\documentclass[onecollarge]{svjour2}       
\smartqed  
\usepackage{amsfonts,latexsym,eucal,amsmath,amssymb} 
\usepackage{graphicx}
\usepackage{psfrag}
\usepackage{color}
\usepackage{amsbsy}
\usepackage{amssymb}       
\usepackage{multicol}
\usepackage{wrapfig}
\usepackage{blindtext}
\usepackage{lineno}
\usepackage{comment}
\usepackage[square,numbers,sort&compress]{natbib} 
\usepackage{braket}
\usepackage{authblk}
\usepackage{fullpage}
\DeclareTextFontCommand{\textmyfont}{\myfont}
\newenvironment{myfont}{\fontfamily{phv}\selectfont}{\par}
\usepackage{array}
\newcolumntype{C}[1]{>{\centering\let\newline\\\arraybackslash\hspace{0pt}}m{#1}}

\usepackage{pifont}
\newcommand{\cmark}{\ding{51}}%
\newcommand{\xmark}{\ding{55}}%

\def\a{{\bf a}}
\def\b{{\bf b}}

\def\dd{\mbox{d}}

\def\q{{\bf q}}

\def\x{{\bf x}}
\def\y{{\bf y}}

\def\sb{{\bf s}}

\def\n{{\bf n}}

\def\HS{\hspace{0.5mm}}

\journalname{Journal}

\begin{document}

\title{A Hierarchical Kinetic Theory of Birth, Death and Fission in Age-Structured 
Interacting Populations}
\titlerunning{Kinetics of age-structured populations}

\author{Tom Chou \and Chris D Greenman}

\institute{
T. Chou \at Depts. of Biomathematics and Mathematics, UCLA, Los Angeles, CA 90095-1766\\
\email{tomchou@ucla.edu}
\and
C. Greenman \at School of Computing Sciences, University of East Anglia, Norwich, UK, NR4 7TJ \\
\email{C.Greenman@uea.ac.uk}
}

\date{Received: date / Accepted: date}

\maketitle

\begin{abstract}
We study mathematical models describing the evolution of stochastic
age-structured populations. After reviewing existing approaches, we
develop a complete kinetic framework for age-structured interacting
populations undergoing birth, death and fission processes in
spatially dependent environments.  We define the full probability
density for the population-size age chart and find results under
specific conditions. Connections with more classical models are also
explicitly derived. In particular, we show that factorial moments for
non-interacting processes are described by a natural generalization of
the McKendrick-von Foerster equation, which describes mean-field
deterministic behavior. Our approach utilizes mixed-type,
multidimensional probability distributions similar to those employed
in the study of gas kinetics and with terms that satisfy BBGKY-like
equation hierarchies.

%
\end{abstract}

\maketitle


\section{Introduction}

Ageing is an important controlling factor in populations with
individuals that range in size from single cells to multicellular
organisms. Age-dependent population dynamics, where birth and death
rates depend on an organism's age, are important in quantitative
models of demography \cite{KEYFITZBOOK}, biofilm formation
\cite{AYATI2007}, stem cell differentiation \cite{SUN2013,ROSHAN2014},
and lymphocyte proliferation and death \cite{ZILMAN2010}. In
cell-based applications, replication is governed by a time-dependent
cell cycle \cite{QU2003,CELLCYCLE2014,AGING0}, while for higher
organisms, the ability to give birth depends on their maturation
time. For applications involving small numbers of individuals, a
stochastic description of the age-structured population is also
desirable.  Finding a practical mathematical framework that captures
age structure, intrinsic stochasticity, and interactions in a
population would be useful for modeling many applications.

Standard frameworks for analyzing age-structured populations include
Leslie matrix models \cite{LESLIE1945,LESLIE1948,MATPOPBOOK}, which
discretizes ages into discrete bins and the continuous-age
McKendrick-von Foerster equation and first studied by McKendrick
\cite{MCKENDRICK,KEYFITZ} and subsequently von Foerster
\cite{VONFOERSTER}, Gurtin and MacCamy \cite{Gurtin1,Gurtin2}, and
others \cite{IANNELLI1995,WEBB2008}.  These approaches describe
deterministic dynamics; stochastic fluctuations in population size are
not incorporated. On the other hand, intrinsic stochasticity and
fluctuations in total population are naturally studied via the
Kolmogorov master equation \cite{VANKAMPEN2011, FPTREVIEW}. However,
the structure of the master equation implicitly assumes exponentially
distributed event (birth and death) times, precluding it from being
used to describe age-dependent rates or age structure within the
population. Evolution of the generating function associated with the
probability distribution for the entire population have also been
developed \cite{BELLMANHARRIS,REID1953,SHONKWILER1980,CHOUJTB}.
Although this approach, the Bellman-Harris equation, allows for
age-dependent event rates, an assumption of independence precludes
population-dependent event rates.
%
%
More recent methods \cite{Jagers2000,Klebaner2014,Hong2011,Hong2013}
have utilized Martingale approaches, which have been used mainly to
investigate the asymptotics of age structure, coalescents, and
estimation of Malthusian growth rate parameters. 

Thus, it is desirable to have a complete mathematical framework that
can resolve the age structure of a population at all time points,
incorporate stochastic fluctuations, and be straightforwardly adapted
to treat nonlinear interactions such as those arising in populations
constrained by a carrying capacity \cite{VERHULST1838,VERHULST1845}.
In a recent publication \cite{Greenman1}, we took a first step in this
direction by formulating a full kinetic equation description that
captures the stochastic evolution of the entire age-structured
population and interactions between individuals.  Here, we
generalize the kinetic equation approach introduced in
\cite{Greenman1} along two main directions.  First, we quantify the
corrections to the mean-field equations by showing that the factorial
moments of the stochastic fluctuations follow an elegant
generalization of the McKendrick-von Foerster equation. Second, we
show how the methods in \cite{Greenman1} can be extended to
incorporate fission processes, where single individuals
instantaneously split into two identical zero-age offspring. These
methods are highlighted with cell division and spatial models. We also
draw attention to the companion paper \cite{Greenman2}, where quantum
field theory techniques developed by Doi and Peliti
\cite{Doi1,Doi2,Peliti} are used to address the same problem,
providing alternative machinery for age-structured modeling.

In the next section, we give a detailed overview of the different
techniques currently employed in age-structured population
modeling. In Section \ref{ANALYSIS}, we use previous results
\cite{Greenman1} to show how the moments of age-structured population
size obey a generalized McKendrick-von Foerster equation. In Section
4, we develop the kinetic theory for branching processes involving
fission. In Section 5, we demonstrate how our theory can be applied to
a microscopic model of cell growth. In Section 6, we demonstrate how
to incorporate spatial effects. Conclusions complete the paper.

\begin{figure}[t]
\begin{center}
\includegraphics[width=11cm]{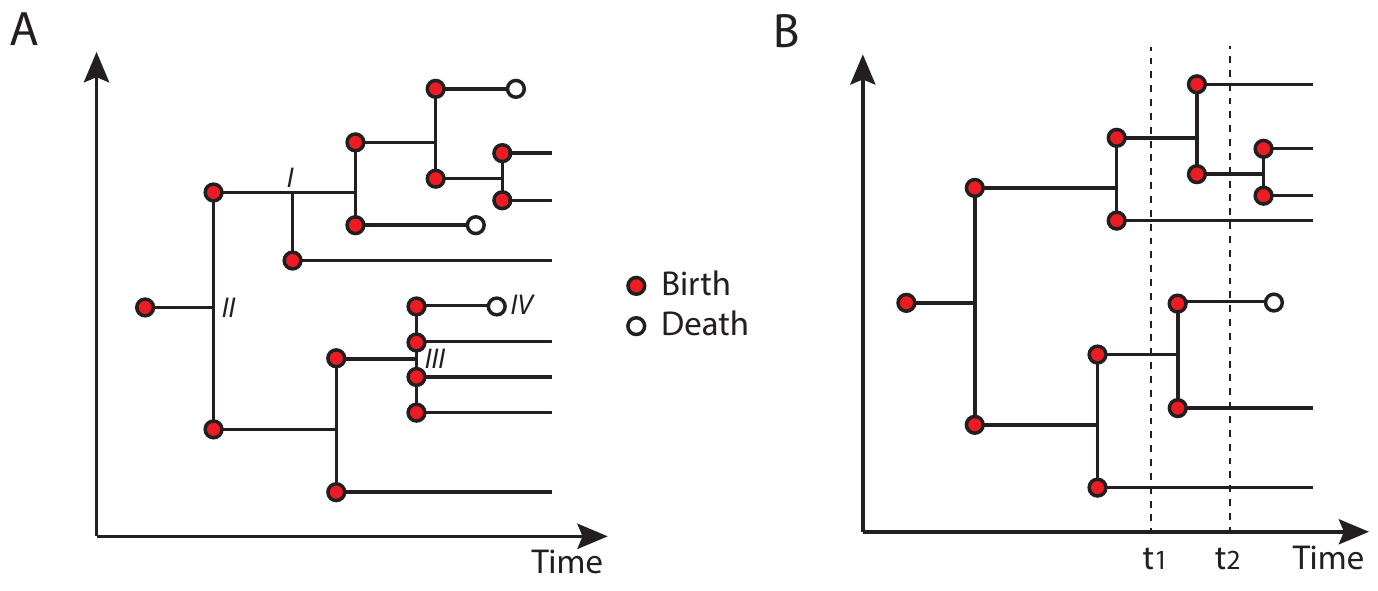}
\caption{\textsf{A}: A general branching process. $\mathsf{I}$ indicates a
  \emph{budding} or \emph{simple} birth process, where the parental
  individual produces a single offspring (a `singlet') without
  death. $\mathsf{II}$ indicates \emph{binary fission}, where a parent
  dies at the same moment two newborn \emph{twins} occur (a
  `doublet'). $\mathsf{III}$ indicates a more general fission event
  with four offspring (a `quadruplet'). $\mathsf{IV}$ indicates death,
  which can be viewed as fission with zero offspring. \textsf{B}: A binary
  fission process such as cell division. At time $t_1$ we have four
  individuals; two sets of twins. At time $t_2$ we have six
  individuals; two pairs of twins and two singlets.}
\label{MicroMod}
\end{center}
\end{figure}


\section{Age-Structured Population Modelling}

Here we review, compare, and contrast existing techniques of
population modeling: the McKendrick-von Foerster equation, the master
equation, the Bellman-Harris equation, Leslie matrices, Martingale
methods, and our recently introduced kinetic approach
\cite{Greenman1}.


\subsection{McKendrick-von Foerster Equation}

It is instructive to first outline the basic structure of the classical
McKendrick-von Foerster deterministic model as it provides a
background for a more complete stochastic picture. First, one defines
$\rho(a,t)$ such that $\rho(a,t)\dd a$ is the expected number of
individuals with an age within the interval $[a,a+\dd a]$. The total
number of organisms at time $t$ is thus $n(t) = \int_{0}^{\infty}
\rho(a,t)\dd a$. Suppose each individual has a rate of giving birth
$\beta(a)$ that is a function of its age $a$. For example, $\beta(a)$
may be a function peaked around the time of \texttt{M} phase in a cell
cycle or around the most fecund period of an organism. Similarly, its
rate of dying $\mu(a)$ can also depend on its age, and will typically 
increase with age.

The McKendrick-von Foerster equation is most straightforwardly derived
by considering the total number of individuals with age in $[0, a]$:
$N(a,t) = \int_{0}^{a}\rho(y,t)\dd y$. The number of births per unit
time from all individuals into the population of individuals with
age in $[0,a]$ is $B(t) = \int_{0}^{\infty} \beta(y)\rho(y,t)\dd y$,
whilst the number of deaths per unit time within this cohort is
$D(a,t) =\int_{0}^{a}\mu(y)\rho(y,t)\dd y$. Within a small time window
$\varepsilon$, the change in $N(a,t)$ is

\begin{equation}
N(a+\varepsilon,t+\varepsilon) - N(a,t) = \int_{t}^{t+\varepsilon}B(s)\dd s
- \int_{0}^{\varepsilon}D(a+s,t+s)\dd s.
\end{equation}
In the $\varepsilon \to 0$ limit, we find

\begin{align}
{\partial N(a,t)\over \partial t} + 
{\partial N(a,t)\over \partial a} & = \int_{0}^{a} \dot{\rho}(y,t)\dd y + \rho(a,t)  = B(t) -\int_{0}^{a} \mu(y)\rho(y,t)\dd y.
\label{INTEGRAL1}
\end{align}
Upon taking ${\partial \over \partial a}$ of Eq.~\ref{INTEGRAL1}, we
obtain the McKendrick-von Foerster equation:

\begin{equation}
{\partial \rho(a,t)\over \partial t} + {\partial \rho(a,t)\over \partial a} = 
-\mu(a)\rho(a,t).
\label{MCKENDRICK0}
\end{equation}
The associated boundary condition arises from setting $a=0$ in
Eq.~\ref{INTEGRAL1}:

\begin{equation}
\rho(a=0,t) =  \int_{0}^{\infty}
\beta(y)\rho(y,t)\dd y \equiv B(t).
\label{MCKENDRICKBC}
\end{equation}
Finally, an initial condition $\rho(a,t=0) = g(a)$
completely specifies the mathematical model.

Note that the term on the right-hand side of Eq.~\ref{MCKENDRICK0} 
depends only on death; the birth rate arises in the boundary condition
(Eq.~\ref{MCKENDRICKBC}) since births give rise to age-zero
individuals. These equations can be formally solved using the method
of characteristics. The solution to Eqs.~\ref{MCKENDRICK0} and
\ref{MCKENDRICKBC} that satisfies a given initial condition is

\begin{equation}
\rho(a,t)=
\begin{cases}
g(a-t)\exp\left[-\int_{a-t}^{a}\mu(s)\dd s\right],
& \text{$a\ge t$}.\\
B(t-a)\exp\left[-\int_{0}^{a}\mu(s)\dd s\right],
& \text{$a<t$}.
\label{MVFSOL}
\end{cases}
\end{equation}

To explicitly identify the solution, we need to calculate the fecundity function
$B(t)$. By substituting the solution in Eq.~\ref{MVFSOL} into
the boundary condition of Eq.~\ref{MCKENDRICKBC} and defining the
propagator $U(a_1,a_2)\equiv\exp\left[-\int_{a_1}^{a_2}\mu(s)\dd s
\right]$, we obtain the following Volterra integral equation:

\begin{equation}
B(t) = \int_0^t B(t-a)U(0,a) \dd a + \int_0^\infty g(a)U(a,a+t) \dd a.
\end{equation}
After Laplace-transforming with respect to time, 
we find 

\begin{equation}
\hat{B}(s) = \hat{B}(s)\mathcal{L}_s\left\{U(0,t)\right\}
+\int_0^\infty g(a)\mathcal{L}_s\left\{U(a,a+t)\right\} \dd a.
\end{equation}
Solving the above for $\hat{B}(s)$ and inverse Laplace-transforming,
we find the explicit expression 

\begin{equation}
B(t) = \mathcal{L}_t^{-1}\left\{ \frac{\int_0^\infty
  g(a)\mathcal{L}_s\left\{U(a,a+t)\right\} \dd
  a}{1-\mathcal{L}_s\left\{U(0,t)\right\}} \right\},
\label{MVFBCSOL}
\end{equation}
which provides the complete solution when used in Eq.~\ref{MVFSOL}.
The McKendrick-von Foerster equation is a deterministic model
describing only the expected age distribution of the population.  If
one integrates Eq.~\ref{MCKENDRICK0} across all ages $0\leq a <
\infty$ and uses the boundary conditions, the rate equation for the
total population is $\dot{n}(t) =
\int_{0}^{\infty}(\beta(a)-\mu(a))\rho(a,t)\dd a$. Thus, $n(t)$ will
in general eventually diverge or vanish depending on the details of
$\beta(a)$ and $\mu(a)$. In the special case $\beta(a) = \mu(a)$, the
population is constant.

What is missing are interactions that stabilize the total population.
Eqs.~\ref{MCKENDRICK0} and \ref{MCKENDRICKBC} assume no higher-order
interactions (such as competition for resources, a carrying capacity,
or mating patterns involving pairs of individuals) within the
populations. Within the McKendrick-von Foerster theory, interactions
are typically heuristically incorporated by assuming that the birth
and death rates ($\beta(a;n(t))$ and $\mu(a;n(t))$) can depend on the
total population $n(t)$ \cite{Gurtin1,Gurtin2,CUSHING1994}. However, as shown in
\cite{Greenman1}, this assumption is an uncontrolled approximation and
inconsistent with a detailed microscopic stochastic model of birth and
death.


\subsection{Master Equation Approach}

A popular way to describe stochastic birth-death processes is through
a function $\rho_{n}(t)$ defining the probability that a population
contains $n$ identical individuals at time $t$. The evolution of this
process can then be described by the standard forward continuous-time
master equation \cite{VANKAMPEN2011,FPTREVIEW}

\begin{equation}
{\partial \rho_{n}(t)\over \partial t} = 
-n\left[\beta_n(t)+\mu_n(t)\right]\rho_n(t)
+(n-1)\beta_{n-1}(t)\rho_{n-1}(t) 
+(n+1)\mu_{n+1}(t)\rho_{n+1}(t),
\label{MASTER0}
\end{equation}
where $\beta_{n}(t)$ and $\mu_{n}(t)$ are the birth and death rates,
per individual, respectively. Each of these rates can be
population-size- and time-dependent.  As such, Eq.~\ref{MASTER0}
explicitly includes the effects of interactions. For example, a
carrying capacity can be implemented into the birth rate through the
following form:
 
\begin{equation}
\beta_{n}(t) = \beta_{0}(t)\left(1-{n\over K(t)}\right).
\label{CAPACITYDEF}
\end{equation}
Here we have allowed both the intrinsic birth rate $\beta_{0}(t)$ and
the carrying capacity $K(t)$ to be functions of
time. Eq.~\ref{MASTER0} can be analytically or numerically solved
via generating function approaches, especially for simple functions
$\beta_{n}$ and $\mu_{n}$.

Since $\rho_{n}(t)$ only describes the total number of individuals at
time $t$, it cannot resolve the distribution of ages within the
fluctuating population. Another shortcoming is the implicit assumption
of exponentially distributed waiting times between birth and death
events. The times since birth of individuals are not tracked. General waiting
time distributions can be incorporated into a master equation approach
by assuming an appropriate number of internal ``hidden'' states, such as
the different phases in a cell division cycle \cite{CELLCYCLE2014}. After all
internal states have been sequentially visited, the system makes a
change to the external population-size state. The waiting time between
population-size changes is then a multiple convolution of the
exponential waiting-time distributions for transitions along each set
of internal states. The resultant convolution can then be used to
approximate an arbitrary waiting-time distribution for the effective
transitions between external states. It is not clear, however, how to
use such an approach to resolve the age structure of the population.


\subsection{Bellman-Harris Fission Process}
\label{BHP}

The Bellman-Harris process \cite{BELLMANHARRIS, REID1953, JAGERS1968,
  SHONKWILER1980, CHOUJTB} describes fission of a particle into any
number of identical daughters, such as events $\mathsf{II}$,
$\mathsf{III}$, and $\mathsf{IV}$ in Fig. \ref{MicroMod}A. Unlike the
master equation approach, the Bellman-Harris branching process
approach allows interfission times to be arbitrarily
distributed. However, it does not model the budding mode of birth indicated by
process $\mathsf{I}$ in Fig.~\ref{MicroMod}A, nor does it capture
interactions (such as carrying capacity effects) within the
population. In such a noninteracting limit, the Bellman-Harris fission
process is most easily analyzed using the generating function $F(z,t)$
associated with the probability $\rho_{n}(t)$, defined as

\begin{equation}
F(z,t) \equiv \sum_{n=0}^{\infty} \rho_{n}(t) z^{n}.
\label{F0}
\end{equation}

We assume an initial condition consisting of a single, newly born
parent particle, $\rho_{n}(0) = \delta_{n,1}$. If we also assume the
first fission or death event occurs at time $\tau$, we can define
$F(z,t\vert \tau)$ as the generating function conditioned on the first
fission or death occurring at time $\tau$ and write $F$ recursively
\cite{LJSALLEN,NEY1972,HARRIS} as:

\begin{equation}
F(z,t\vert \tau) = \begin{cases}
z, & t<\tau, \\
A(F(z,t-\tau)), & t \ge \tau,
\end{cases}
\hspace{10mm}
A(x) = \sum_{m=0}^{\infty} a_{m} x^{m}.
\label{BH0}
\end{equation}
The function $A$ encapsulates the probability $a_{m}$ that a particle
splits into $m$ identical particles upon fission, for each
non-negative integer $m$. For binary fission, we have $A(x) = (1-a_{2})
+ a_{2}x^{2}$ since $\sum_{m=0}^{\infty}a_{m} = 1$. Since this overall
process is semi-Markov \cite{HONG}, each daughter behaves as a new
parent that issues its own progeny in a manner statistically equivalent to and
independent from the original parent, giving rise to the
compositional form in Eq.~\ref{BH0}. We now weight $F(z,t\vert \tau)$
over a general distribution of waiting times between splitting events,
$g(\tau)$, to find

\begin{align}
F(z,t) & \displaystyle \equiv \int_{0}^{\infty}F(z,t \vert \tau)g(\tau)\dd \tau \nonumber\\
& = z \int_{t}^{\infty}\!\!\!g(\tau)\dd \tau + 
\int_{0}^{t}\!A(F(z,t-\tau))g(\tau)\dd \tau.
\label{BELLMANHARRIS}
\end{align}

The Bellman-Harris branching process \cite{NEY1972, FOK_PROOF} is thus
defined by two parameter functions: $a_{m}$, the vector of progeny
number probabilities, and $g(\tau)$, the probability density function
for waiting times between branching events for each particle. The
probabilities $\rho_n(t)$ can be recovered using a contour integral
surrounding or a Taylor expansion about the origin:

\begin{equation}
\rho_{n}(t) =\frac{1}{2\pi i}\oint_{C}\frac{F(z,t)}{z^{k+1}} \dd z
= {1\over n!}{\partial^{n} F(z,t)\over \partial z^{n}} \bigg|_{z=0}.
\end{equation}

Note that Eq.~\ref{BELLMANHARRIS} allows one to incorporate an
arbitrary waiting-time distribution between events, a feature that is
difficult to implement in the master equation (Eq.~\ref{MASTER0}). An
advantage of the branching process approach is the ease with which
general waiting-time distributions, multiple species, and immigration
can be incorporated. However, it is limited in that an independent
particle assumption was used to derive Eq.~\ref{BELLMANHARRIS}, where
the statistical properties of the entire process starting from one
parent were assumed to be equivalent to those started by each of the
identical daughters born at time $\tau$. This assumption of
independence precludes treatment of interactions within the
population, such as those giving rise to carrying capacity.
%
%
More importantly, the Bellman-Harris equation is expressed purely in
terms of the generating function for the total population size and
cannot resolve age structure within the population.


\subsection{Leslie Matrices}

Leslie matrices \cite{LESLIE1945,LESLIE1948} have been used to resolve
the age structure in population models \cite{SUN2013,REID1953,
  THESIS1998,DIFFUSION2009,GETZ1984,COHEN,LESLIE1945,LESLIE1948,CUSHING}. These
methods essentially divide age into discrete bins and are implemented
by assuming fixed birth and death rates within each age bin. Such
approaches have been applied to models of stochastic harvesting
\cite{GETZ1984,COHEN} and fluctuating environments
\cite{LANDE1988,LANDE2005} and are based on the following linear
construction, iterated over a single time step:

\begin{equation}
\begin{bmatrix} n_0\\n_1\\ \vdots \\ n_{N-1}  \end{bmatrix}_{t+1}=
\begin{bmatrix}
f_0 & f_1 & \hdots & f_{N-2} & f_{N-1} \\ 
s_0 & 0 & \hdots & 0 & 0 \\
0 & s_1 & \hdots & 0 & 0 \\
\vdots & \vdots & \ddots & \vdots & \vdots\\
0 & 0 & \hdots & s_{N-2} & 0
\end{bmatrix}
\cdot
\begin{bmatrix} n_0\\n_1\\ \vdots \\ n_{N-1}  \end{bmatrix}_{t}.
\end{equation}
The value $n_i$ indicates the population size in age group $i$; $f_i$
is the mean number of offspring arriving to age group $0$ from a
parent in age group $i$; and $s_i$ is the fraction of individuals
surviving from age group $i$ to $i+1$. These models have the advantage
of being based upon algebraic linearity, which enables many features
of interest to be investigated analytically
\cite{MATPOPBOOK}. However, they are inherently deterministic
(although they can be used to study extrinsic environmental noise) and
the discretization of such models results in an approximation. Thus, a
fully continuous stochastic model is desirable.


\subsection{Martingale Approaches}

Relatively recent investigations have used Martingale approaches to
model age-structured stochastic processes. These methods stem from
stochastic differential equations and Dynkin's formula \cite{OKSENDAL}
and considers general processes of the form $F(f(\a_n(t)))$, where the
vector $\a_n(t)$ represents the time dependent age-chart of the
population with variable size $n$; $f$ is a symmetric function of the
individual ages; and $F$ is a generic function of interest. A
Martingale decomposition of the following form results,

\begin{equation}
F(f(\a_n(t))) = F(f(\a_n(0))) + \int_0^t {\cal G} F(f(\a_n;s)) \dd s + M_t^{(f,F)},
\end{equation}
where the operator ${\cal G}$ captures the mean behavior, and the
stochastic behavior is encoded in  the local Martingale process
$M_t^{(f,F)}$ \cite{Jagers2000}. Such analyses have enabled several
features of general birth-death processes, including both budding and
fission forms of birth to be quantified. Specifically, the Malthusian
growth parameter can be explicitly determined, along with the
asymptotic behavior of the age-structure. More recently there have
been results related to coalescents and extinction of these processes
\cite{Klebaner2014,Hong2011,Hong2013}. However, we will show the
utility of obtaining the probability density of the entire age chart
of the population which allows efficient computations in transient
regimes. The kinetic approach first developed in \cite{Greenman1}
introduces machinery to accomplish this.


   
\subsection{Kinetic Theory}

A brief introduction to the current formulation of our kinetic theory
approach to age-structured populations can be found in
\cite{Greenman1}.  The starting point is a derivation of a
variable-dimension Liouville equation for the complete probability
density function $\rho_n(\a_n;t)$ describing a stochastic,
interacting, age-structured population subject to simple birth and
death.  Variables in the theory include the population size $n$, time
$t$, and the vector $\a_n=(a_1,a_2,\hdots,a_n)$ representing the
complete age chart for the $n$ individuals. If we randomly label the
individuals $1,2,\hdots,n$, then $\rho_n(\a_n;t)\dd \a_n$ represents
the probability that the $i^\mathrm{th}$ individual has age in the
interval $[a_i,a_i+\dd a_i]$. Since individuals are considered
indistinguishable, $\rho_n(\a_n;t)$ is invariant under any permutation
of the age-chart vector $\a_{n}$. These functions are analogous to
those used in kinetic theories of gases \cite{MCQUARRIE}. Their
analysis in the context of age-structured populations builds on the
Boltzmann kinetic theory of Zanette \cite{ZANETTE} and results in a
BBGKY-like (Bogoliubov-Born-Green-Kirkwood-Yvon) hierarchy of
equations:

\begin{align}
\displaystyle {\partial \rho_{n}(\a_{n};t)\over \partial t} + 
& \sum_{j=1}^{n}{\partial \rho_{n}(\a_{n};t)\over \partial a_{j}} = 
-\rho_{n}(\a_{n};t)\sum_{i=1}^{n} \gamma_{n}(a_{i}) 
 + (n+1)\!\int_{0}^{\infty}\!\!\mu_{n+1}(y)
\rho_{n+1}(\a_{n},y;t)\dd y,
\label{RHO0} 
\end{align}
where $\gamma_n(a) = \beta_n(a) +\mu_n(a)$ and the age variables are
separated from the time variable by the semicolon. The associated
boundary condition is given by

\begin{equation}
\begin{array}{l}
n \rho_{n}(\a_{n-1},0;t) = \rho_{n-1}(\a_{n-1};t)\beta_{n-1}(\a_{n-1}).
\label{BCRHO}
\end{array}
\end{equation}

\begin{table}[t]
\centering
\caption{Advantages and disadvantages of different frameworks for
  stochastic age-structured populations. `Stochastic' indicates that
  the model resolves probabilities of configurations of the population
  `Age-dependent rates' indicates whether or not a model takes into
  account birth, death, or fission rates that depend on an individuals
  age (time after its birth). `Age-structured Populations' indicates
  whether or not the theory outputs the age structure of the ensemble
  population. `Age Chart Resolved' indicates whether or not a theory
  outputs the age distribution of all the individuals in the
  population. `Interactions' indicates whether or not the approach can
  incorporate population-dependent dynamics such as that arising from
  a carrying capacity, or from birth processes involving multiple
  parents.  `Budding' and `Fission' describes the model of birth and
  indicates whether the parent lives or dies after
  birth. \textsuperscript{1}Birth and death rates in the McKendrick-von Foerster 
equation can be made explicit functions of the total populations size, which must be
self-consistently solved \cite{Gurtin1,Gurtin2}.
\textsuperscript{2}Leslie matrices discretize age groups and
  are an approximate method. \textsuperscript{3}Martingale methods do
  not resolve the age structure explicitly, but utilize 
  rigorous machinery. \textsuperscript{4}The kinetic approach
  for fission is addressed later in this work, but not in \cite{Greenman1}.}
\label{TabTech}
\begin{tabular}{|m{2.3cm}||C{1.6cm}|C{1.4cm}|C{1.4cm}|C{1.4cm}|C{1.3cm}|C{1.3cm}|C{1.1cm}|}
\hline Theory & Stochastic & Age-dependent rates & Age-structured
Populations & Age Chart Resolved & Interactions & Budding & Fission
\\ \hline\hline\vspace{1mm} Verhulst Eq. & \xmark & \xmark & \xmark &
\xmark & \cmark & \xmark & \xmark \\ \hline\vspace{1mm} McKendrick Eq.
& \xmark & \cmark & \cmark & \xmark & \cmark & \cmark\textsuperscript{1}& \xmark
\\ \hline\vspace{1mm} Master Eq. & \cmark & \xmark & \xmark & \xmark &
\cmark & \cmark & \cmark \\ \hline\vspace{1mm} Bellman-Harris & \cmark
& \cmark & \xmark & \xmark& \xmark & \xmark & \cmark
\\ \hline\vspace{1mm} Leslie Matrices& \xmark &
\cmark\textsuperscript{2} & \cmark & \xmark & \cmark & \xmark & \xmark
\\ \hline\vspace{1mm} Martingale & \cmark & \cmark &
\xmark\textsuperscript{3} & \xmark & \cmark & \cmark & \cmark
\\ \hline\vspace{1mm} Kinetic Theory & \cmark & \cmark & \cmark &
\cmark & \cmark & \cmark & \cmark\textsuperscript{4} \\ \hline
\end{tabular}
\end{table}

Note that because $\rho_n(\a_{n-1},0;t)$ is symmetric in the age
arguments, the zero can be placed equivalently in any of the $n$ age
coordinates. The birth rate function is quite general and can take
forms such as $\beta_{n-1}(\a_{n-1})=
\sum_{i=1}^{n-1}\beta_{n-1}(a_i)$ for a simple birth process or
$\sum_{1 \le i < j \le n-1}\beta_{n-1}(a_i,a_j)$ to represent births
arising from interactions between pairs of individuals.

Equation \ref{RHO0} applies only to the budding or simple
mode of birth such as event $\mathsf{I}$ in Fig.~\ref{MicroMod}A. In
\cite{Greenman1} we derived analytic solutions for
$\rho_{n}(\a_{n};t)$ in pure death and pure birth processes and showed
that when the birth and death rates are constant, the hierarchy of
equations reduces to a single master equation.
Characterizing all the remaining higher moments of the distribution
remains an outstanding problem. Moreover, methods to tackle
fission modes of birth such as those shown in Fig.~\ref{MicroMod}B 
were not developed. These are the two main areas of focus of the
present work. Before analyzing these problems, we summarize the pros
and cons of the different techniques in Table \ref{TabTech}.


\section{Analysis of Simple Birth-Death Processes}
\label{ANALYSIS}

We now revisit the simple budding birth and death process and extend
the the kinetic framework introduced in \cite{Greenman1}.  We first
show that the factorial moments for the density $\rho_n(\a_n;t)$
satisfy a generalized McKendrick-von Foerster equation. We also
explicitly solve  Eqs.~\ref{RHO0} and \ref{BCRHO}, and 
derive for the first time an exact general solution for
$\rho_n(\a_n;t)$.


\subsection{Moment Equations}
\label{MOMENTEQNS}

The McKendrick-von Foerster equation has been shown to correspond to a
mean-field theory of age-structured populations \cite{Greenman1}.
Specifically, if we construct the expected population density $X(a,t) = \sum_n
n\rho_n^{(1)}(a;t)$, we obtain the McKendrick-von Foerster equation
for $X(a,t)$. This leaves open the problem of determining the
age-structured variance (and higher-order moments) of the population
size.

In \cite{Greenman1}, we derived the marginal $k-$dimensional
distribution functions defined by integrating $\rho_n(\a_n;t)$ over
$n-k$ age variables:

\begin{equation}
\rho_{n}^{(k)}(\a_{k};t) \equiv \int_{0}^{\infty}\!\dd a_{k+1}
\ldots\int_{0}^{\infty}\!\dd a_{n}\, \rho_{n}(\a_{n};t).
\end{equation}
The symmetry properties of $\rho_{n}(\a_{n};t)$ indicate that it is
immaterial which of the $n-k$ age variables are integrated out. 
From Eq.~\ref{RHO0}, we then obtained

\begin{align}
\displaystyle {\partial \rho_{n}^{(k)}(\a_{k};t) \over \partial t}  + 
\sum_{i=1}^{k}{\partial \rho_{n}^{(k)} (\a_{k};t)\over \partial a_{i}} 
= &  \displaystyle - \rho_{n}^{(k)}(\a_{k};t)\sum_{i=1}^{k}\gamma_{n}(a_{i}) \nonumber \\
\: & + \left({n-k \over n}\right)\rho_{n-1}^{(k)}(\a_{k};t)
\sum_{i=1}^{k}\beta_{n-1}(a_{i}) \nonumber \\
\: &  +{(n-k)(n-k-1)\over n}\int_{0}^{\infty}\beta_{n-1}(y)
\rho_{n-1}^{(k+1)}(\a_{k},y;t)\dd y \nonumber \\
\: & -(n-k)\int_{0}^{\infty}\gamma_{n}(y)
\rho_{n}^{(k+1)}(\a_{k},y;t)\dd y \label{RHOK1} \\ 
\: & \displaystyle + (n+1)\int_{0}^{\infty}\mu_{n+1}(y)
\rho_{n+1}^{(k+1)}(\a_{k},y;t)\dd y.\nonumber
\label{RHOK1}
\end{align}
Similarly, integrating the boundary condition in Eq.~\ref{BCRHO} over
$n-k$ of the (nonzero) variables, gives, for simple birth processes
where $\beta_n(\a_m)= \sum_{i=1}^m \beta_n(a_i)$,

\begin{equation}
\rho_n^{(k)}(\a_{k-1},0;t) =
\frac{1}{n}\rho_{n-1}^{(k-1)}(\a_{k-1};t)\sum_{i=1}^{k-1}
\beta_{n-1}(a_i)+\frac{n-k}{n}\int_0^\infty
\rho_{n-1}^{(k)}(\a_{k-1},y;t)\beta_{n-1}(y)\dd
y.
\label{BCMARG}
\end{equation}
We now show how to use these marginal density
equation hierarchies and boundary conditions to derive an equation for
the $k^\mathrm{th}$ moment of the age density.

For $k=1$, $\rho_{n}^{(1)}(a;t)\dd a$ is the probability that we have
$n$ individuals and that if one of them is randomly chosen, it will
have age in $[a,a+\dd a]$. Therefore, the probability that we have $n$
individuals, and there exists an individual with age in $[a,a+\dd a]$,
is $n\rho_{n}^{(1)}(a; t)\dd a$. Summing over all possible population
sizes $n\geq 1$ yields the probability $\rho(a,t)\dd a =
\sum_nn\rho_n^{(1)}(a;t)\dd a$ that the system contains an individual
with age in the interval $[a,a+\dd a]$. More generally,
$n^k\rho_n^{(k)}(\a_k;t)\dd \a_k$ is the probability that there are
$n$ individuals, $k$ of which can be labelled such that the
$i^\mathrm{th}$ has age within the interval $[a_i,a_i+\dd
  a_i]$. Summing over the possibilities $n \ge k$, we thus introduce
factorial moments $X^{(k)}(\a_k;t)$ and moment functions
$Y^{(k)}(\a_k;t)$ as:

\begin{align}
X^{(k)}(\a_k;t) \equiv \sum_{n=k}^{\infty}(n)_k\rho_{n}^{(k)}(\a_k;t)
\equiv \sum_{\ell=0}^k s(k,\ell)Y^{(\ell)}(\a_\ell;t),
\nonumber\\
Y^{(k)}(\a_k;t) \equiv \sum_{n=k}^{\infty}{n^k}\rho_{n}^{(k)}(\a_k;t)
\equiv \sum_{\ell=0}^k S(k,\ell)X^{(\ell)}(\a_\ell;t).
\label{XANDY}
\end{align}
Here $(n)_k = n(n-1)\hdots(n-(k-1))=k!{n\choose k}$ is the Pochhammer
symbol, and $s(k,\ell)$, $S(k,\ell)$ are Stirling numbers of the first
and second kind \cite{STANLEY2012,STANLEY2001}. Although we are
primarily interested in the functions $Y^{(k)}(\a_k;t)$, the factorial
moments $X^{(k)}(\a_k;t)$ will prove to be analytically more
tractable. One can then easily interchange between the two moment
types by using the polynomial relationships involving Stirling
numbers.

After multiplying Eq.~\ref{RHOK1}
by $(n)_k$ and summing over all $n \ge k$, we find

\begin{align}
\frac{\partial X^{(k)}}{\partial t}+ 
\sum_{i=1}^{k}\frac{\partial X^{(k)}}{\partial a_i}  
+\sum_{n \ge k}(n)_k\rho_{n}^{(k)} & \sum_{i=1}^{k}\gamma_{n}(a_{i}) = 
\sum_{n-1 \ge k}(n-1)_k\rho_{n-1}^{(k)}
\sum_{i=1}^{k}\beta_{n-1}(a_{i}) \nonumber \\
&  +\int_{0}^{\infty}\sum_{n-1 \ge k+1}(n-1)_{k+1}\rho_{n-1}^{(k+1)}(\a_{k},y;t)\beta_{n-1}(y)\dd y
\nonumber \\
&  -\int_{0}^{\infty}\sum_{n \ge k+1}(n)_{k+1}\rho_{n}^{(k+1)}(\a_{k},y;t)\gamma_{n}(y)\dd y
\nonumber \\
&  +\int_{0}^{\infty}\sum_{n+1 \ge k+1}(n+1)_{k+1}\rho_{n+1}^{(k+1)}(\a_{k},y;t)\mu_{n+1}(y)\dd y,
\label{MVFGC}
\end{align}
where, for simplicity of notation, the arguments $(\a_k;t)$ have been
suppressed from $\rho_{n}^{(k)}$ and $X^{(k)}$.  In the case where the
birth and death rates $\beta_n(a)=\beta(a)$ and $\mu_n(a)=\mu(a)$ are
independent of the sample size, significant cancellation occurs and we
obtain the elegant equation

\begin{equation}
\frac{\partial X^{(k)}}{\partial t}+ 
\sum_{i=1}^{k}\frac{\partial X^{(k)}}{\partial a_i}  
+X^{(k)}\sum_{i=1}^{k}\mu(a_{i}) = 0.
\label{MVFG}
\end{equation}
When $k=1$, one recovers the classical McKendrick-von Foerster
equation describing the mean-field behavior after stochastic
fluctuations are averaged out.  Equation \ref{MVFG} is a natural
generalization of the McKendrick-von Foerster equation and provides
all the age-structured moments arising from the population size
fluctuations.  If the birth and death rates, $\beta_{n}$ and
$\mu_{n}$, depend on the population size, one has to analyze the
complicated hierarchy given in Eq.~\ref{MVFGC}.

To find the boundary conditions associated with
Eq.~\ref{MVFG}, we combine the definition of $X^{(k)}$ with the
boundary condition in Eq.~\ref{BCMARG} and obtain

\begin{align}
X^{(k)}(\a_{k-1},0;t) = & \sum_{n \ge k}(n)_k\rho_n^{(k)}(\a_{k-1},0;t)\nonumber\\
= & X^{(k-1)}(\a_{k-1};t)\beta(\a_{k-1})+\int_0^\infty X^{(k)}(\a_{k-1},y;t)\beta(y) \dd y.
\label{FMBC}
\end{align}
Setting $X^{(0)}\equiv 0$, we recover the boundary
condition associated with the classical McKendrick-von Foerster
equation. For higher-order factorial moments, the full solution to the
$(k-1)^\mathrm{st}$ factorial moment $X^{(k-1)}(\a_{k-1};t)$ is required
for the boundary condition to the $k^\mathrm{th}$ moment
$X^{(k)}(\a_{k-1},0;t)$.

Specifically, consider the second factorial moments and assume the
solution $X^{(1)}\equiv Y^{(1)}$ to the McKendrick-von Foerster
equation is available (from \textit{e.g.}, Eq.~\ref{MVFSOL}).  In the
small interval $\dd a$, the term $Y^{(1)}\dd a$ is the Bernoulli
variable for an individual having an age in the interval
$[a,a+da]$. Thus, in an extended age window $\Omega$, we heuristically
obtain the expectation

\begin{equation}
\mathrm{E}(Y_\Omega(t)) = \sum_{da \in \Omega}Y_{da}(t)=\int_\Omega Y^{(1)}(a;t) \dd a,
\label{MNEQ}
\end{equation}
where $Y_\Omega(t)$ is the stochastic random variable describing the
number of individuals with an age in $\Omega$ at time $t$.
Using an analogous argument for the variance, we find

\begin{equation}
\mathrm{Var}(Y_\Omega(t))=\sum_{da,db \in \Omega}
Cov(Y_{da},Y_{db})=\int_{\Omega^2}Y^{(2)}(a,b;t)\dd a \dd
b-\int_\Omega Y^{(1)}(a;t)\dd a \cdot \int_\Omega Y^{(1)}(b;t) \dd b.
\label{VAREQ}
\end{equation}
Thus, the second moment $Y^{(2)}$ allows us to describe the variation
of the population size within any age region of interest. Similar
results apply for higher order correlations. We focus then on deriving
a solution to $Y^{(2)}$ and determining the variance of
population-size-age-structured random variables.  Eq.~\ref{MVFG} for
general $k$ is readily solved using the method of
characteristics

\begin{equation}
X^{(k)}(\a_k;t)=X^{(k)}(\a_k-m;t-m)\prod_{j=1}^kU(a_j-m,a_j),
\label{FMOMSOL}
\end{equation}
where the propagator is defined as $U(a,b) \equiv \exp\left[-\int_a^b
  \mu(s) \dd s\right]$. We can now specify $X^{(k)}$ in terms of
boundary conditions and initial conditions by selecting $m =
\min\left\{\a_k,t\right\}$.  Since $X^{(k)}(\a_k;t)\equiv
X^{(k)}(\pi(\a_k);t)$ is invariant to any permutation $\pi$ of its age
arguments, we have only two conditions to consider.  The initial
condition $X^{(k)}(\a_k;0)=g(\a_k)$ encodes the initial correlations
between the ages of the founder individuals and is assumed to be
given. From Eq.~\ref{XANDY}, $X^{(k)}(\a_k;0)$ must be a symmetric
function in the age arguments.  A boundary condition of the form
$X^{(k)}(\a_{k-1},0;t)\equiv B(\a_{k-1};t)$ describes the fecundity of
the population through time. This is not given but can be determined
in much the same way that Eq.~\ref{MVFBCSOL} was derived.

To be specific, consider a simple pure birth (Yule-Furry) process
($\beta(a) = \beta$, $\mu(a) = 0$) started by a single individual. The
probability distribution of the initial age of the parent individual
is assumed to be exponentially distributed with mean $\lambda$. Upon
using transform methods similar to those used to derive
Eq.~\ref{MVFBCSOL}, we obtain the following factorial moments (see
Appendix A for more details):

\begin{equation}
X^{(1)}(a;t) = 
\begin{cases}
\lambda e^{-\lambda(a-t)}, & t<a\\
\beta e^{\beta (t-a)}, & t>a
\end{cases},
\hspace{10mm}
X^{(2)}(a,b;t) = 
\begin{cases}
0, & t<a<b\\
\lambda\beta e^{-\lambda(b-a)}e^{(\lambda+\beta)(t-a)}, & a<t<b\\
2\beta^2 e^{-\beta(b-a)}e^{2\beta(t-a)}, & a<b<t
\end{cases}.
\label{APPASTUFF}
\end{equation}
We have given $X^{(2)}(a,b;t)$ for only $a<b$
since the region $a>b$ can be found by imposing 
symmetry of the age arguments in $X^{(2)}$.
After using Eq.~\ref{XANDY} to convert $X^{(1)}$
and $X^{(2)}$ into $Y^{(1)}$ and $Y^{(2)}$, we can use Eqs.~\ref{MNEQ}
and \ref{VAREQ} to find age-structured moments, particularly the
mean and variance for the number of individuals that have age in the
interval $[a,b]$:

\begin{align}
\mathrm{E}(Y_{[a,b]}(t)) & = 
\begin{cases}
e^{\lambda(t-a)}-e^{\lambda(t-b)}, & t<a<b\\
e^{\beta(t-a)}-e^{\lambda(t-b)}, & a<t<b\\
e^{\beta(t-a)}-e^{\beta(t-b)}, & a<b<t,
\end{cases}\\
\mathrm{Var}(Y_{[a,b]}(t)) & = 
\begin{cases}
e^{2\lambda t}(e^{-\lambda a}-e^{-\lambda b})(-e^{-\lambda a}+e^{-\lambda b}+e^{-\lambda t}), & t<a<b\\
(e^{\beta(t-a)}-e^{\lambda(t-b)}) (e^{\beta(t-a)} +e^{\lambda(t-b)}-1), & a<t<b\\
e^{2\beta t}(e^{-\beta a}-e^{-\beta b})(e^{-\beta a}-e^{-\beta b}+e^{-\beta t}), & a<b<t.
\end{cases}
\end{align}
Note that in the limits $a \to 0$ and $b\to \infty$, we recover the
expected exponential growth of the total population size
$\mathrm{E}(Y_{[0,\infty]}) = e^{\beta t}$ for a Yule-Furry
process. We also recover the known total population
variance $\mathrm{Var}(Y_{[0,\infty]}) = e^{\beta t}(e^{\beta t}-1)$.


\subsection{Full Solution}

Equation \ref{RHO0} defines a set of coupled linear
integro-differential equations in terms of the density
$\rho_n(\a_n;t)$. In \cite{Greenman1}, we derived explicit analytic
expressions for $\rho_n(\a_n;t)$ in the limits of pure death and pure
birth.  Here, we outline the derivation of a formal expression for the
full solution.  It will prove useful to revert to the following
representation for the density:

\begin{equation}
f_n(\a_n;t) \equiv n!\rho_n(\a_n;t).
\label{REP}
\end{equation}
If $\a_n$ is restricted to the ordered region such that $a_1
\le a_2 \le \hdots \le a_n$, $f_n$ can be interpreted as the probability
density for age-ordered individuals (see \cite{Greenman1} for more
details). We will consider $f_n$ as a distribution over
$\mathbb{R}^n$; however, its total integral ($n!$) is not unity
and it is not a probability density. We can use Eq.~\ref{REP} to
switch between the two representations, but simpler analytic
expressions for solutions to Eq.~\ref{RHO0} result when $f_n(\a_n;t)$
is used.

\begin{figure}[t!]
  \centering
  \setlength{\unitlength}{0.1\textwidth}
  \begin{picture}(12,4)
    \put(0.5,0){\includegraphics[width=7.5cm]{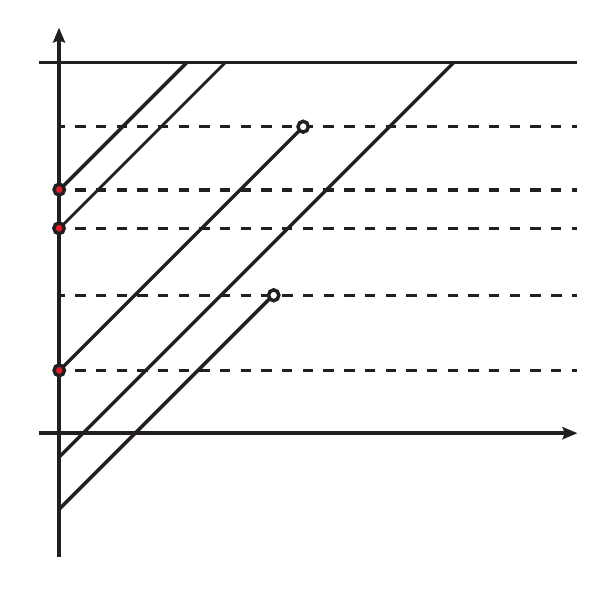}}
    \put(5,1.00){Age}
    \put(0.8,4.70){Time}
    \put(0.6,3.72){$s_1$}
    \put(0.6,3.17){$b_2$}
    \put(0.6,2.83){$b_1$}
    \put(0.6,2.27){$s'_1$}
    \put(0.6,1.7){$y_1$}
    \put(0.6,0.95){$b'_1$}
    \put(0.6,0.57){$y'_1$}
    \put(0.53,4.35){$t$}
    \put(0.5,1.20){$0$}
    \put(5.6,4.2){$\hat{f}_3(\b_2,b'_1;t)=\hat{U}(\b_2,b'_1;s_1,t)\cdot$}
    \put(5.6,3.72){$\displaystyle\int_{b_2}^{t} \!\!\dd s_1 \int_0^{s'_1} \!\!\! \dd y_1 
\hat{U}(\b_2,y_1,b'_1;b_2,s_1)\mu_4(s_1-y_1)\cdot$}
    \put(5.6,3.17){$\hat{U}(b_1,y_1,b'_1;b_1,b_2) \beta_3(b_2-[b_1,y_1,b'_1])\cdot$}
    \put(5.6,2.83){$\hat{U}(y_1,b'_1;s'_1,b_1) \beta_2(b_1-[y_1,b'_1])\cdot$}
    \put(5.6,2.27){$\displaystyle\int_{y_1}^{b_1} \!\!\dd s'_1 \int_{-\infty}^0 \!\!\! 
\dd y'_1 \hat{U}(y_1,b'_1,y'_1;y_1,s'_1)\mu_3(s'_1-y'_1)\cdot$}
    \put(5.6,1.66){$\hat{U}(b'_1,y'_1;0,y_1)\beta_2(y_1-[b'_1,y'_1])\cdot$}
    \put(5.6,1.20){$\hat{f}_2(b'_1,y'_1;0)$}
  \end{picture}
  \caption{A sample birth death process over the time interval
    $[0,t]$.  Red and white circles indicate births and deaths within
    $[0,t]$. The variables $b_i>0$ and $b'_j<0$ denote
    TOBs of individuals present at time $t$, while $y_i>0, y'_j<0$,
    and $s_i,s'_j \in [0,t]$ indicate birth and death times of
    individuals who have died by time $t$. Terms
    arising from application of the recursion in Eq.~\ref{RECC} and
    boundary condition of Eq.~\ref{BCRHO} are given to the right.}
\label{TIMEAGE}
\end{figure}
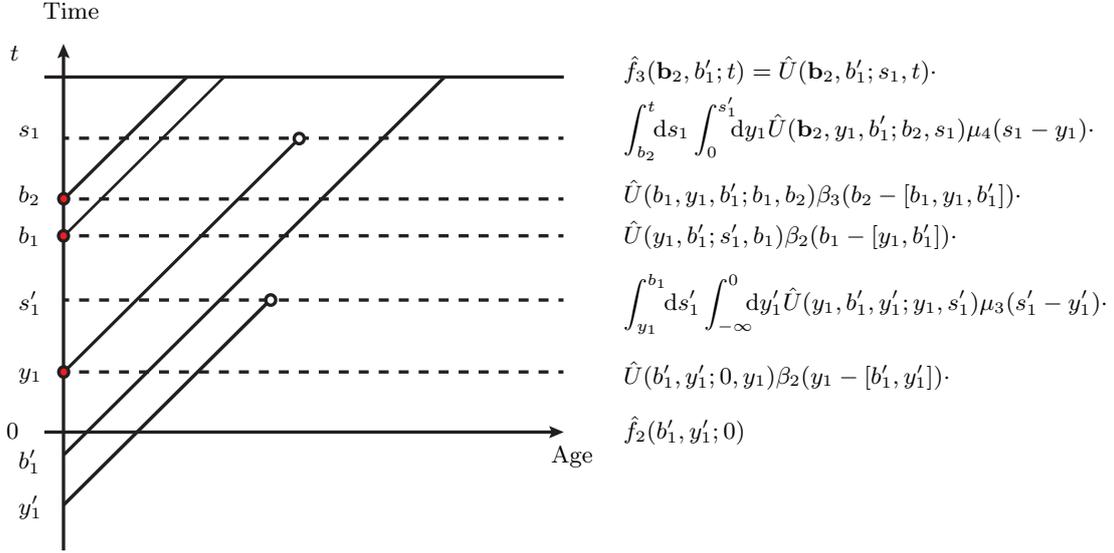

To find general solutions for $f_n(\a_n;t)$ expressed in terms of an
initial distribution $f_{m}(\cdot;0)$, we replace $\rho_n(\a_n;t)$
with $f_n(\a_n;t)/n!$ in Eq.~\ref{RHO0} and use the method of
characteristics to find a solution. Examples of characteristics are
the diagonal timelines portrayed in Fig.~\ref{TIMEAGE}. So far,
everything has been expressed in terms of the natural parameters of
the system; the age $\a_{n}$ of the individuals at time $t$. However,
$\a_{n}$ varies in time complicating the analytic expressions. If we
index each characteristic by the time of birth (TOB) $b=t-a$ instead
of age $a$, then $b$ is fixed for any point $(a,t)$ lying on a
characteristic, resulting in further analytic simplicity. We use the
following identity to interchange between TOB and age representations:

\begin{equation}
\hat{f}_n(\b_n;t) \equiv f_n(\a_n;t),
\hspace{10mm}\b_n = t - \a_n.
\label{AANDB}
\end{equation}
We will abuse notation throughout our derivation by
identifying $t-\a_n\equiv [t-a_1,t-a_2,\hdots,t-a_n]$.  The method of
characteristics then solves Eq.~\ref{RHO0} to give a solution of the
following form, for any $t_0 \ge \max \{\b_n\}$
\begin{equation}
\hat{f}_n(\b_n;t)= \hat{f}_n(\b_n;t_0)\hat{U}_n(\b_n;t_0,t)
+\displaystyle\int_{t_0}^t  \dd s  \int_{-\infty}^s \hspace{-3mm}\dd y \HS 
\hat{U}_n(\b_n;s,t) 
\hat{f}_{n+1}(\b_n,y;s)\mu_{n+1}(s-y).
\label{RECC}
\end{equation}
This equation is defined in terms of a propagator
$\hat{U}_n(\b_m;t_0;t) \equiv U_n(\a_m;t_0;t)$ that represents the
survival probability over the time interval $[t_0,t]$, for $m$
individuals born at times $\b_m$, in a population of size $n$,

\begin{equation}
\hat{U}_n(\b_m;t_0,t)=\exp\left[-\sum_{i=1}^{m} \int_{t_0}^t\gamma_{n}(s-b_i)\dd s\right],
\end{equation}
where we have again used the definition $\gamma_n(a)=\beta_n(a)+\mu_n(a)$.
The propagator $\hat{U}$ satisfies certain translational properties:

\begin{align}
\hat{U}_n(\b_m;t_0,t) & =\prod_{i=1}^m\hat{U}_n(b_i;t_0,t),
\label{TRANS1}\\
\hat{U}_n(\b_m;t_0,t) & =\hat{U}_n(\b_m;t_0,t')\cdot \hat{U}_n(\b_m;t',t).
\label{TRANS2}
\end{align}
The solution $\hat{f}_n$ applies to any region of phase space where
$t_0\ge\max(\b_n)$. If $t_0=\max(\b_n)$, say $t_0=b_n$, then
we must invoke the boundary conditions of Eq.~\ref{BCRHO} to replace
$\hat{f}_n (\b_{n-1},b_n;b_n)$ with $\hat{f}_{n-1}
(\b_{n-1};b_n)\beta_{n-1}(b_n-\b_{n-1})$, where we have used and will henceforth use
the notation 

\begin{align}
\beta_{n-1}(b_n-\b_{n-1}) & \equiv  \beta_{n-1}(b_n-[b_{1}, b_{2},\ldots,b_{n-1}])\nonumber \\
\: & \equiv \sum_{i=1}^{n-1}\beta_{n-1}(b_n-b_i).
\end{align}
Eq.~\ref{RECC} is then used to propagate $\hat{f}_{n-1}
(\b_{n-1};b_n)$ backwards in time. To obtain a general solution, we
need to repeatedly back-substitute Eq.~\ref{RECC} and the associated
boundary condition, resulting in an infinite series of
integrals. However, each term in the resultant sum can be represented
by a realization of the birth-death process. We represent any such
realization across time period $[0,t]$, such as that given in
Fig.~\ref{TIMEAGE}, as follows.

Let $\b_m \in [0,t]$ and $\b'_n<0$ denote the TOBs for $m$ individuals
born in the time interval $[0,t]$, and $n$ founder individuals, all
alive at time $t$.  Next, define $\y_k \in [0,t]$ and $\y'_\ell<0$ to be the
TOBs of $k$ individuals born in the time interval $[0,t]$ and $\ell$
founder individuals, respectively. Here, all $k+\ell$ individuals are
assumed to die in the time window $[0,t]$. Their corresponding times
of death are defined as $\sb_k$ and $\sb'_\ell$, respectively. Thus, there will
be $n+\ell$ individuals alive initially at time $t=0$ and $m+n$
individuals alive at the end of the interval $[0,t]$.

Next, consider the realization in Fig.~\ref{TIMEAGE}, where we start
with the two individuals at time $0$ with TOBs $b'_1$ and $y'_1$. The
individual with TOB $b'_1$ survives until time $t$, while the individual with
TOB $y'_1$ dies at time $s'_1$. Within the time interval $[0,t]$ there
are three more births with TOBs $b_1$, $b_2$ and $y_1$, the last of
which has a corresponding death time of $s_1$, resulting in three
individuals in total that exist at time $t$. 

To express the distribution $\hat{f}_3(\b_2,b'_1;t)$ in terms of the
initial distribution $\hat{f}_2(b'_1,y'_1;0)$, conditional upon
three birth and two death events ordered such that
$0<y_1<s'_1<b_1<b_2<s_1<t$, we start with the distribution
$\hat{f}_2(b'_1,y'_1;0)$. Just prior to the first birth time $y_1$,
we have two individuals, so that $\hat{f}_3(\cdot;y_1^{-})\equiv 0$
and Eq.~\ref{RECC} yields $\hat{f}_2(b'_1,y'_1;y_1^{-}) =
\hat{f}_2(b'_1,y'_1;0) \hat{U}(b'_1,y'_1;0,y_1)$ (the death
term does not contribute).  To describe a birth at time $y_1$, we use
the boundary condition of Eq.~\ref{BCRHO} to construct
$\hat{f}_3(b'_1,y'_1,y_1;y_1)=
\hat{f}_2(b'_1,y'_1;y_1^{-})\beta_2(y_1-[b'_1,y'_1])$. 

Immediately after $y_1$ and before the next death occurs at time
$s'_1$, three individuals exist and $\hat{f}_2(\cdot;y_1^{+})\equiv
0$. Now, only the death term in Eq.~\ref{RECC} contributes and
\begin{equation}
\hat{f}_2(y_1,b'_1;b_1^{-}) = \int_{y_1}^{b_1}\dd s'_1
\int_{-\infty}^0 \dd y'_1 \\ \hat{U}(y_1,b'_1,y'_1;y_1,s'_1)
\mu_3(s'_1-y'_1) \hat{f}_3(y_1,b'_1,y'_1;s'_1).
\end{equation}
Continuing this counting, we find the product of terms displayed on
the right-hand side of Fig.~\ref{TIMEAGE}. Next, we use the
translational properties indicated in Eqs.~\ref{TRANS1} and
\ref{TRANS2} to combine the propagators associated with
Fig.~\ref{TIMEAGE} into one term:
$\hat{U}(y'_1;0,s'_1)\hat{U}(b'_1;0,t)
\hat{U}(y_1;y_1,s_1)\hat{U}(b_1;b_1,t)\hat{U}(b_2;b_2,t)$.  In other
words, each birth-death pair $(y,s)$ is propagated along the time
interval it survives; from $\max\{y,0\}$ to $\min\{s,t\}$. For
example, the individual with TOB $b'_1<0$ survives across the entire
timespan $[0,t]$, whereas the individual with TOB $y_1$ is born and
dies at times $y_1$ and $s_1$. These two individuals are propagated by
the terms $U(b'_1;0,t)$ and $U(y_1;y_1,s_1)$, respectively.  Provided
the order $0<y_1<s'_1<b_1<b_2<s_1<t$ is preserved and the values
$b'_1,y'_1<0$ are negative, the form of the integral expressions in
Fig.~\ref{TIMEAGE} are preserved.

\begin{figure}
\begin{center}
\hspace{-2mm}\includegraphics[width=16.2cm]{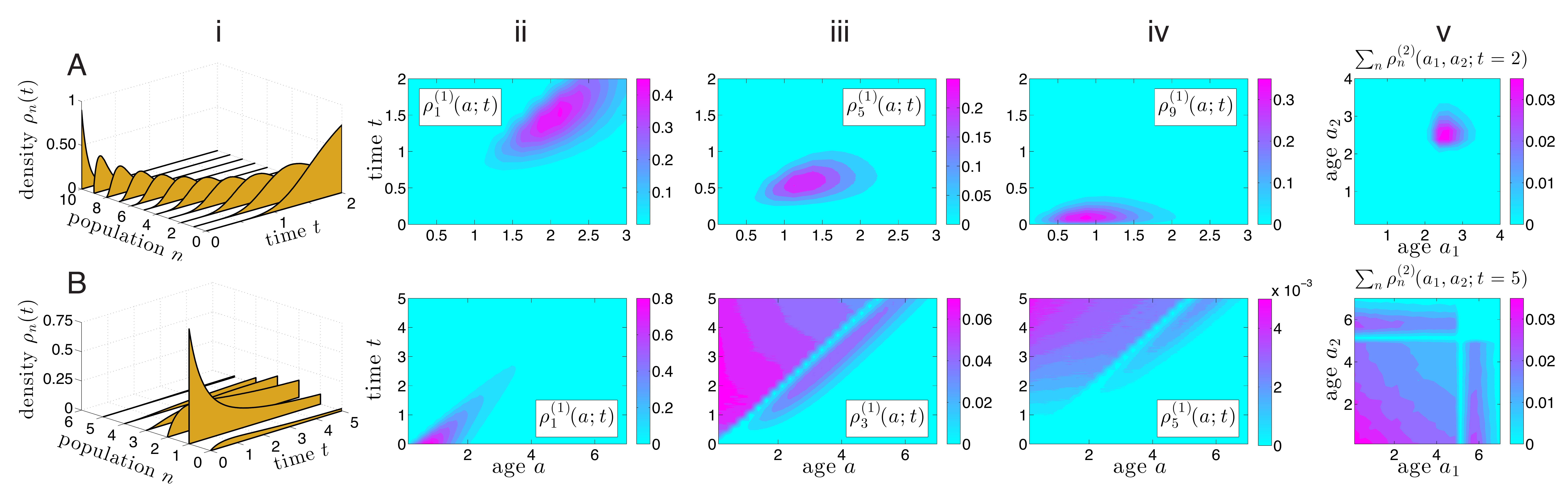}
\vspace{-1mm}
\caption{Monte-Carlo simulations of densities in age- and
  capacity-dependent birth-death processes. Row \textsf{A} are results
  for a death-only process with a linear death rate function $\mu(a) =
  a$.  We initiated all simulations from $10$ individuals with initial age
  drawn from distribution $P(a) = 128a^{3}e^{-4a}/3$.
  In row \textsf{B}, we consider the a budding-only birth process with
  a carrying capacity $K=5$ (in Eq.~\ref{CAPACITYDEF}). Here,
  simulations were initiated with a single parent individual with an
  initial age also drawn from the distribution $P(a)$. In (\textsf{i}),
  we plot the total number density
  $\rho_{n}^{(0)}(t)= \int \dd \a \HS \rho_n(\a;t)$ for both
  processes. We also plot the single-particle density function
  $\rho_{n=1,5,9}(a;t=2)$ for the pure death process in
  \textsf{A}(\textsf{ii-iv}) and $\rho_{n=1,3,5}(a;t=5)$ for the
  limited budding process in \textsf{B}(\textsf{ii-iv}). Finally, the
  population-summed two-point correlations functions
  $\sum_{n}\rho_{n}^{(2)}(a_{1},a_{2};t)$ for pure death and pure
  budding are shown in panels \textsf{A}(\textsf{v}) and
  \textsf{B}(\textsf{v}).
\label{MONTECARLO}}
%
\end{center}
\end{figure}

After summing across all realizations $C_{m,k,\ell}$ (the
configuration in Fig.~\ref{TIMEAGE} is one member of $C_{2,1,1}$) of
the possible orderings of the birth and death times $\b_m$, $\y_k$,
$\y'_\ell$, $\sb_k$ and $\sb'_\ell$, we can write the general solution
to Eq.~\ref{RECC} in the form

\begin{align}
\hat{f}_{m+n}(\b_m,\b'_n;t) = &
\sum_{k,\ell=0}^\infty\sum_{C_{m,k,\ell}} 
\int_{-\infty}^0\dd \y'_\ell\cdot
\int_{t^{-}(\y_k)}^{t^{+}(\y_k)}\dd \y_k\cdot
\int_{t^{-}(\sb_k)}^{t^{+}(\sb_k)}\dd \sb_k\cdot
\int_{t^{-}(\sb'_\ell)}^{t^{+}(\sb'_\ell)}\dd \sb'_\ell\cdot
\hat{f}_{n+l}(\b'_n,\y'_\ell;0)\cdot
\nonumber\\
& \prod_{i=1}^mU(b_i;b_i,t)\cdot 
\prod_{i=1}^kU(y_i;y_i,s_i)\cdot
\prod_{i=1}^nU(b'_i;0,t)\cdot
\prod_{i=1}^\ell U(y'_i;0,s'_i)\prod_{i=1}^m \beta_{N(b_i)}(b_i-A(b_i))\cdot
\nonumber\\
&
\prod_{i=1}^k \beta_{N(y_i)}(y_i-A(y_i))\cdot
\prod_{i=1}^k \mu_{N(y_i)}(s_i-y_i)\cdot
\prod_{i=1}^\ell \mu_{N(y'_i)}(s'_i-y'_i).
\label{GENSOL}
\end{align}
The terms $t^{-}(\x)$ and $t^{+}(\x)$ refer to the times below and
above $\x$ relative to the ordering of times $\b_m$, $\y_k$,
$y'_\ell$, $\sb_k$ and $\sb'_k$. For example, in Fig.~\ref{TIMEAGE},
$t^{-}(\b_2)=[s'_1,b_1]$ and $t^{+}(\b_2)=[b_2,s_1]$ represent the
lower and upper bounds of the vector $\b_2=[b_1,b_2]$ found from the
ordering $0<y_1<s'_1<b_1<b_2<s_1$. The term $A(x)$ represents the
vector of TOBs of the individuals alive just prior to time $x$. The
term $N(x)$ represents the number of individuals alive just prior to
time $x$. 

Although analytic and complete, the solution given in Eq.~\ref{GENSOL}
is unwieldy and difficult to implement.  One can truncate the sum to
remove low probability contributions, such as realizations containing
improbable numbers of intermediary births and deaths, and perform
numerical integration. However, this approach also rapidly becomes
infeasible as the dimensions increase.  Therefore, we explore the
general solution via event-based Monte-Carlo simulation. We initialize
the process with a number of samples obtained from an initial
distribution. Each sample is represented by a vector $\b_n$ of birth
times and is propagated forward in time.  A timestep is chosen to be
sufficiently small such that at most one birth or death event occurs
within it, after which the vector $\b_n$ is updated.  This process is
continued until the required time has been reached.  Although the high
dimensionality makes it difficult to sample enough realizations to
sufficiently explore the distribution $f_n(\a_n;t)$, lower dimensional
marginal distributions such as $f_n^{(0)}(\cdot;t)$,
$f_n^{(1)}(a_1;t)$ and $f_n^{(2)}(a_1,a_2;t)$, and their counterparts
$\rho_{n}$, can be sufficiently sampled.

Figures~\ref{MONTECARLO}\textsf{A} and \textsf{B} show results from
simulations of a pure death and a pure birth process, respectively.
In Fig.~\ref{MONTECARLO}\textsf{A} we assumed a population-independent
linear death rate $\mu(a) = a$ and initiated the pure death process
with 10 individuals with initial ages drawn from a gamma distribution with unit mean 
and standard deviation $\frac{1}{2}$.  Fig.~\ref{MONTECARLO}\textsf{A(i)}
shows the simulated density which decreases in $n$ with time.
Figs.~\ref{MONTECARLO}\textsf{A(ii-iv)} show that the weight of the
reduced single-particle density function shifts to longer times and
higher ages as the system size $n$ is decreased.  The sum over the
population of the symmetric two-point correlation
$\rho_{n}^{(2)}(a_{1}, a_{2};t=2)$ is shown in
Fig.~\ref{MONTECARLO}\textsf{A(v)}. The observed structure indicates
no correlations in the death only process and the peak at $a_{1}=a_{2}
\approx 2.6$ reflects the fact that older individuals die faster,
shifting the mean age slightly below the initial age plus the elapsed
time (1+2=3). Results from Monte-Carlo simulations of a pure birth
process with growth rate $\beta_{0} = 1$ and carrying capacity $K=5$
(Eq.~\ref{CAPACITYDEF}). Here, we initiated the simulations with one
individual with age drawn from the same gamma distribution $P(a) = 128
a^{3}e^{-4a}/3$.  In this case, the reduced single-particle density
exhibits peaks arising from both from the initial distribution and
from birth (Fig.~\ref{MONTECARLO}\textsf{B(ii-iv)}). The two-point
correlation function $\sum_{n=0}^{\infty}\rho_{n}^{(2)}(a_{1},
a_{2};t=5)$ exhibits a similar multimodal structure as shown in
\textsf{(v)}.

In all simulations at least 400,000 trajectories were aggregated and
the results are in good agreement with analytic solutions to
Eq.~\ref{RHO0}. Similar analytic results can be obtained using
Doi-Peliti second quantization methods, as is demonstrated in the
companion paper \cite{Greenman2}. In particular, the age-structured
population-size function $\rho_n(t)$ is expanded into a similar sum,
where each term can be interpreted two ways: as an element in a
perturbative expansion and also represented as a Feynman diagram in a
path integral expansion. The moment equations from Section
\ref{MOMENTEQNS} that generalize the McKendrick equation can also be
derived using second quantization.


\section{Age-Structured Fission-Death Processes}

We now derive a kinetic theory for a binary fission-death process, as
depicted in Fig.~\ref{MicroMod}B. We find a hierarchy of kinetic
equations, analogous to Eqs.~\ref{RHO0} and \ref{BCRHO}, and determine
the mean behavior.

\subsection{Extended Liouville Equation for Fission-Death}

The binary fission-death process is equivalent to a birth-death
process except that parents are instantaneously replaced by
\textit{two} newborns. The process can also be thought of as a budding
process in which the parent is instantaneously renewed.  In order to
describe both twinless individuals (singlets) and twins (a doublet),
we have to double the dimensionality of our density functions. For
example, in Fig.~\ref{MicroMod}B at time $t_1$, we have two pairs of
distinct twins, with four individuals having two ages, whereas at time
$t_2$ we have two singlets and two doublets. Thus, we define the ages
of current singlets and twins by $\a_m$ and $\a_n'$, respectively,
where $m$ is the number of singlets and $n$ the number of pairs of
twins.  Transforming to the time-of-birth (TOB) representation, we
define the TOB of current singlets and twins as $\x_m=t-\a_m$ and
$\y_n=t-\a_n'$, respectively. For simplicity, we will assume that no
simple birth processes occur and that particles grow in number only
through fission. The function $\beta_{m,n}(a)$ is defined as the
age-dependent fission rate of an individual (whether a singlet or a
doublet) of age $a$ when the system contains $m$ singlets and $n$
doublets. Similarly, we have death rate $\mu_{m,n}(a)$, and event rate
$\gamma_{m,n}(a) = \beta_{m,n}(a)+ \mu_{m,n}(a)$. We suppose, for the
moment, that the TOBs are ordered so that $x_1 \le x_2 \le \hdots \le
x_m$ and $y_1 \le y_2 \le \hdots \le y_m$.  The quantity
$f_{m,n}(\x_{m};\y_{n})\dd \x_m \dd \y_n$ is then the probability of
$m$ singlets with ordered TOBs in $[\x_m,\x_m+\dd \x_m]$ and $n$ twin
pairs with ordered TOBs in $[\y_n,\y_n+\dd \y_n]$.  The density
$f_{m,n}$ satisfies the following equation:

\begin{equation}
\begin{array}{l}
\displaystyle {\partial f_{m,n}(\x_{m};\y_{n};t)\over \partial t}
+f_{m,n}(\x_{m};\y_{n};t)\left[\sum_{i=1}^{m}\gamma_{m,n}(t-x_{i})+2\sum_{j=1}^{n}\gamma_{m,n}(t-y_{j})\right]
= \\[11pt] \: \hspace{2cm} \displaystyle
\sum_{i=0}^{m}\int_{x_{i}}^{x_{i+1}}
f_{m+1,n}(\x_{i},z,\x_{i+1,m};\y_{n};t)\mu_{m+1,n}(t-z)\dd z \\[11pt]
\: \hspace{2cm} + \displaystyle
2\sum_{i=1}^{m}f_{m-1,n+1}(\x_m^{(-i)};\y_{i},x_{i},
\y_{i+1,n};t)\mu_{m-1,n+1}(t-x_{i}),
\label{FISSIONF0}
\end{array}
\end{equation}
where the partial age vectors are defined as $\x_{i,j} =
(x_i,\hdots,x_j)$ and the singlet age vector, doublet age vector, and
time arguments are separated by semicolons. The term
$\x_m^{(-i)}=(x_1,\hdots,x_{i-1},x_{i+1},\hdots,x_m)$ represents the
vector of all $m$ singlet TOBs, except for the $i^\mathrm{th}$ one. 
The first term on the right hand side of Eq.~\ref{FISSIONF0}
represents the death of a singlet particle with an unknown TOB $z$ in
the interval $[x_i,x_{i+1}]$, while the second term describes the
death of any one of two individuals in a pair of twins (with TOB
$x_i$).

The associated boundary conditions are:

\begin{align} 
f_{m,n}(\x_{m-1},t; \y_{n};t) & = 0,\\
f_{m,n}(\x_{m}; \y_{n-1},t;t) & = 2\sum_{i=1}^{m}
f_{m-1,n}(\x_m^{(-i)};\y_{n-1},x_{i};t)\beta_{m-1,n}(t-x_{i}) \nonumber\\
& \hspace{3mm}+\sum_{i=0}^{m} \int_{x_{i}}^{x_{i+1}}f_{m+1,n-1}(\x_{i},z,\x_{i+1,m};
\y_{n};t)\beta_{m+1,n-1}(t-z)\dd z.
\end{align}
The first term on the right-hand side above represents the fission of
one of a pair of twins, generating a new pair of twins of age zero
(TOB $t$), and leaving behind a singlet with TOB $x_{i}$. The second
term represents the fission (and removal) of a singlet with unknown
TOB $z$, giving rise to an additional pair of twins of age zero.

We now let $\x_m$ and $\y_n$ be unordered TOB vectors, and extend
$f_{m,n}$ to the domain $\mathbb{R}^{m+n}$ by defining
$f_{m,n}(\x_{m};\y_{n};t) =
f_{m,n}(\mathcal{T}(\x_{m});\mathcal{T}(\y_{n});t)$, where $\mathcal{T}$
is the ordering operator. Note that $f_{m,n}$ is not a probability
distribution under this extension; however,
$\rho_{m,n}(\x_{m};\y_{n};t)\dd \x_m \dd \y_n = {1 \over m!
  n!}f_{m,n}(\x_{m};\y_{n};t)\dd \x_m \dd \y_n$ can be interpreted as
the probability that we have a population of $m$ singlets and $n$
pairs of twins, such that if we randomly label the singlets
$1,2,\hdots,m$ and the doublets $1,2,\hdots,n$, the $i^\mathrm{th}$
singlet has age in $[x_i,x_i+\dd x_i]$ and the $j^\mathrm{th}$ doublet
have age in $[x_j,x_j+\dd x_j]$.  The density $\rho_{m,n}$ obeys

\begin{align}
\frac{\partial \rho_{m,n}(\x_{m};\y_{n};t)}{\partial t}
+\rho_{m,n}(\x_{m};\y_{n};t)& \left[\sum_{i=1}^{m}\gamma_{m,n}(t-x_{i})+2\sum_{j=1}^{n}\gamma_{m,n}(t-y_{j})\right] = 
\nonumber\\
& (m+1)\int_{-\infty}^{t}
\rho_{m+1,n}(\x_{m},z;\y_{n};t)\mu_{m+1,n}(t-z)\dd z 
\nonumber\\
& + 2\left({n+1\over m}\right)\sum_{i=1}^{m}
\rho_{m-1,n+1}(\x_m^{(-i)};\y_{n},x_{i};t)\mu_{m-1,n+1}(t-x_{i}),
\label{FISSRHO}
\end{align}
with associated boundary condition

\begin{align}
\rho_{m,n}(\x_{m-1},t; \y_n;t) & = 0,\nonumber\\
\rho_{m,n}(\x_{m}; \y_{n-1},t;t) & = {2\over
  m}\sum_{i=1}^{m}\rho_{m-1,n}(\x_m^{(-i)};\y_{n-1},x_{i};t)
\beta_{m-1,n}(t-x_{i}) \nonumber\\ 
& \hspace{3mm} +
\left({m+1\over n}\right)\int_{-\infty}^{t}
\rho_{m+1,n-1}(\x_{m},z;\y_{n-1};t) \beta_{m+1,n-1}(t-z)\dd z.
\label{FISSBCRHO}
\end{align}
Equations \ref{FISSRHO} and \ref{FISSBCRHO} provide a complete
probabilistic description of the population of singlets and doublets
undergoing fission and death. We draw attention to the parallel paper
\cite{Greenman2}, where we derive an equivalent hierarchy using
methods used in quantum field theory developed by Doi and Peliti
\cite{Doi1,Doi2,Peliti}.



\subsection{Mean-Field Behavior}

Here, we analyze the mean-field behavior of the fission-death process
by first integrating out unwanted variables from the full density
$\rho_{m,n}(\x_m;\y_n;t)$ to construct marginal or ``reduced''
densities. Successive integrals over any number of the variables
$\x_m$ and $\y_n$ can be performed, giving:

\begin{equation}
\rho_{m,n}^{(k,\ell)}(\x_{k};\y_{\ell};t) \equiv 
\int_{-\infty}^{t}\dd \x'_{m-k}\int_{-\infty}^{t}\dd \y'_{n-\ell}
\rho_{m,n}(\x_k,\x'_{m-k};\y_{\ell},\y'_{n-\ell};t).
\end{equation}

For example, $\rho_{m,n}^{(0,0)}(;;t)$ is the probability of finding
at time $t$, $m$ singlets and $n$ doublets, regardless of
age. After integrating
Eq.~\ref{FISSRHO} we find the double hierarchy of equations

\begin{align}
& \frac{\partial \rho_{m,n}^{(k,\ell)}(\x_{k};\y_{\ell};t)}{\partial t} 
+\rho_{m,n}^{(k,\ell)}(\x_{k};\y_{\ell};t)
\left[ \sum_{i=1}^k\gamma_{m,n}(t-x_i)+ 2\sum_{i=1}^k\gamma_{m,n}(t-x_i) \right]\nonumber\\
& \hspace{10mm}
+ (m-k)\int \rho_{m,n}^{(k+1,\ell)}(\x_{k},z;\y_{\ell};t)
\gamma_{m,n}(t-z) \dd z
+2(n-\ell)\int \rho_{m,n}^{(k,\ell+1)}(\x_{k};\y_{\ell},z;t)
\gamma_{m,n}(t-z) \dd z \nonumber\\
& = (m+1)\int_{-\infty}^t \rho_{m+1,n}^{(k+1,\ell)}(\x_{k},z;\y_{\ell};t) 
\mu_{m+1,n}(t-z) \dd z \nonumber\\
& \hspace{10mm}
+ 2 \left( \frac{n+1}{m} \right)\sum_{i=1}^k
\rho_{m-1,n+1}^{(k-1,\ell+1)}(\x_{k}^{(-i)};\y_{\ell},x_i;t)
\mu_{m-1,n+1}(t-x_i) \nonumber \\
& \hspace{10mm}
+ 2 \left( \frac{n+1}{m} \right)(m-k)\int_{-\infty}^t
\rho_{m-1,n+1}^{(k,\ell+1)}(\x_{k};\y_{\ell},z;t)
\mu_{m-1,n+1}(t-z)  \dd z.
\label{FISSMARG}
\end{align}
Similarly, integrating Eq.~\ref{FISSBCRHO}
yields boundary conditions for the marginal densities:

\begin{align}
\rho_{m,n}^{(k,\ell)}(\x_{k-1},t;\y_{\ell};t) = & 0,
\nonumber\\
\rho_{m,n}^{(k,\ell)}(\x_k;\y_{\ell-1},t;t) = &
\frac{2}{m}\sum_{i=1}^k \rho_{m-1,n}^{(k-1,\ell)}(\x_k^{(-i)},t;\y_{\ell-1},x_i;t) 
\beta_{m-1,n}(t-x_i)
\nonumber\\
& + 2\left(\frac{m-k}{m}\right)\int_{-\infty}^t
\rho_{m-1,n}^{(k,\ell)}(\x_k;\y_{\ell-1},z;t) \beta_{m-1,n}(t-z)\dd z \nonumber\\
& + \left(\frac{m+1}{n}\right)\int_{-\infty}^t
\rho_{m+1,n-1}^{(k+1,\ell-1)}(\x_k,z;\y_{\ell-1};t) \beta_{m+1,n-1}(t-z)\dd z.
\label{FISSMARGBC}
\end{align}

We can now analyze the densities $X(x,t)$ and $Y(y,t)$, where
$X(x,t)\dd x$ is the probability that there exists at time $t$ a
singlet with TOB in $[x,x+\dd x]$ and $Y(y,t) \dd y$ is the
probability that at time $t$ we have one doublet with TOB in
$[y,y+\dd y]$. Analogous to Eq.~\ref{XANDY}, we define

\begin{align}
X(x,t) &  \equiv \sum_{m,n=0}^{\infty}m\rho_{m,n}^{(1,0)}(x;;t) 
= \sum_{m,n=0}^{\infty}m\int_{-\infty}^{t}\dd \x_{m-1}
\int_{-\infty}^{t}\dd \y_{n}\rho_{m,n}(\x_{m-1},x;\y_{n};t), \nonumber \\
Y(y,t) &  \equiv \sum_{m,n=0}^{\infty}n\rho_{m,n}^{(0,1)}(;y;t) =  \sum_{m,n=0}^{\infty}n\int_{-\infty}^{t}\dd \x_{m}
\int_{-\infty}^{t}\dd \y_{n-1} \rho_{m,n}(\x_{m};\y_{n-1},y;t).
\label{FISSMN}
\end{align}
Upon setting $(k,\ell)=(1,0)$ and $(k,\ell)=(1,0)$, we multiply by and sum
over $m$ and $n$ in Eq.~\ref{FISSMARG}, respectively.  If the fission
and death rates $\beta_{m,n}(a)$ and $\mu_{m,n}(a)$ depend on
population size, the resultant expressions are complex hierarchies
which will be difficult to analyze. However, if
$\beta_{m,n}(a)=\beta(a)$ and $\mu_{m,n}(a)=\mu(a)$ are
size-independent, many cancellations occur and the resulting equations
for $X$ and $Y$ simplify significantly

\begin{equation}
\frac{\partial X}{\partial t} =  (2Y-X)\gamma(t-x),
\hspace{10mm}
\frac{\partial Y}{\partial t} = -2Y\gamma(t-x).
\label{MEANFIELDEQNS}
\end{equation}
Similarly, repeating the operation on the boundary conditions in
Eq.~\ref{FISSMARGBC}, we find boundary conditions for $X$ and $Y$:

\begin{equation}
X(t,t) =  0,
\hspace{10mm}
Y(t,t) =  \int_{-\infty}^{t}(X(z,t)+2Y(z,t)) \gamma(t-z)\dd z \equiv B(t).
\label{MEANFIELDBCS}
\end{equation}
Note that if $T=X+2Y$ is the total population density,
Eqs.~\ref{MEANFIELDEQNS} and \ref{MEANFIELDBCS} reduce to
McKendrick-von Foerster-like equations:

\begin{equation}
\frac{\partial T}{\partial t} = -\gamma(t-z)T,
\hspace{10mm}
T(t,t) =  \int_{-\infty}^{t}T(z,t)\gamma(t-z)\dd z.
\end{equation}

To solve Eqs.~\ref{MEANFIELDEQNS} and \ref{MEANFIELDBCS}, we first define
\begin{equation}
U(x;t_1,t_2)=\exp\left[-\int_{t_1}^{t_2}\gamma(s-x) ds\right],
\label{UEQ}
\end{equation}
and find solutions of the form

\begin{align}
X(x,t) & = X(x,t_0)U(x;t_0,t)+2Y(x,t_0)U(x;t_0,t)(1-U(x;t_0,t)),\nonumber\\
Y(x,t) & = Y(x,t_0)U^{2}(x;t_0,t),
\end{align}
provided $t_{0} \geq x$. For an initial time of 
$t=0$, we find, upon setting $t_0=\max\{0,x\}$, 

\begin{equation}
  X(x,t)=\begin{cases}
    2B(x)U(x;x,t)(1-U(x;x,t)), & x>0,\\
    X(x,0)U(x;0,t)+2Y(x,0)U(x;0,t)(1-U(x;0,t)), & x<0,
  \end{cases}
\label{XXSOL}
\end{equation}

\begin{equation}
  Y(x,t)=\begin{cases}
    B(x)U^{2}(x;x,t), & x>0,\\
    Y(x,0)U^{2}(x;0,t), & x<0.
  \end{cases}
\label{YYSOL}
\end{equation}
We now substitute Eqs.~\ref{XXSOL} and \ref{YYSOL} 
into Eqs.~\ref{MEANFIELDBCS} to find a Volterra equation for $B$:

\begin{equation}
B(t) = 2\int_0^{t}B(x)U(x;x,t)
\beta(t-x) \dd x + \int_{-\infty}^0[X(x,0)+2Y(x,0)]U(x;0,t)\beta(t-x) \dd x.
\label{PHIEQN}
\end{equation}
Equation \ref{PHIEQN} along with Eqs.~\ref{XXSOL} and \ref{YYSOL}
constitute a complete solution for the mean density of singlets and
doublets. Eqs.~\ref{XXSOL} and \ref{YYSOL} also show that the total
population density, $T(x,t)= X(x,t)+2Y(x,t)$, takes on a simple form
in terms of $B(t)$:

\begin{equation}
  T(x,t)=\begin{cases}
    2B(t)U(x;x,t), & x>0,\\
    T(x,0)U(x;0,t), & x<0,
  \end{cases}
\label{TFULLFORMULA}
\end{equation}
while the total mean population $T(t)= \int_{0}^{\infty}T(x,t)\dd x$ is given by

\begin{equation}
T(t) = 2\int_0^{t}B(x)U(x;x,t) \dd x + \int_{-\infty}^0 T(x,0)U(x;0,t) \dd x.
\label{MEQN}
\end{equation}
Before analyzing a specific model of the fission-death process, we
will first establish the equivalence of our noninteracting kinetic
theory with the Bellman-Harris fission process (discussed in
Subsection \ref{BHP}) in the mean-field limit.


\subsection{Mean-field Equivalence to the Bellman-Harris Process}

Consider a Bellman-Harris fission process with a branching time
distribution $g(a)$ and a cumulative density function $G(a) =
\int_{0}^{a}g(a')\dd a'$, along with the progeny distribution function
$A(\cdot)$ given in Eq.~\ref{BH0}. The Bellman-Harris model 
in Eq.~\ref{BELLMANHARRIS} can be written in the form

\begin{equation}
F(z,t) = z(1-G(t)) + \int_0^t A(F(z,\tau)) g(t-\tau)\dd \tau.
\end{equation}
If we restrict ourselves to a binary fission process, the progeny
distribution function takes the form $A(y)=a_0+a_2y^2$, where $a_0$
and $a_2 = 1-a_{0}$ are the death and binary fission probabilities,
conditional on an event taking place. Thus, the mean population
defined as
\begin{equation}
T(t) \equiv \left.\frac{\partial F}{\partial z}\right |_{z=1} =
\int_t^{\infty}g(\tau)\dd \tau +
2a_2\int_0^{t}g(t-\tau)T(\tau)\dd \tau
\end{equation}
has the Laplace-transformed solution

\begin{equation}
\tilde{T}(s) = \frac{1}{s}\frac{1-\tilde{g}(s)}{1-2a_2\tilde{g}(s)}.
\label{MLAPLACE}
\end{equation}

We now show that the same result arises from our full noninteracting
(population-independent $\beta(a)$ and $\mu(a)$) kinetic
approach. Since the fission and death rates can be expressed as
$\beta(y)=\frac{a_{2}g(y)}{1-G(y)}$ and $\mu(y)=\frac{a_{0} g(y)}{1-G(y)}$,
Eq.~\ref{UEQ} reduces to $U(x;x,t)=1-G(t-x)$ and
$U(0;0,t)=1-G(t)$. Starting from a single individual at age zero,
Eq.~\ref{MEQN} can be written as

\begin{equation}
T(t) = 2\int_0^{t}B(x)(1-G(t-x)) \dd x + (1-G(t)),
\end{equation}
which has the Laplace-transformed solution

\begin{equation}
\tilde{T}(s) = (2\tilde{B}(s)+1)\frac{1-\tilde{g}}{s}.
\label{TILDETS}
\end{equation}

Similarly, Eq.~\ref{PHIEQN} becomes

%

\begin{equation}
B(t) = 2\int_0^{t}B(x)a_2 g(t-x) \dd x + a_2 g(t),
\end{equation}
with Laplace-transformed solution
\begin{equation}
\tilde{B}(s) = \frac{a_2\tilde{g}(s)}{1-2a_2\tilde{g}(s)}.
\label{PHILAPLACE}
\end{equation}
Substituting Eq.~\ref{PHILAPLACE} in Eq.~\ref{TILDETS} results in
Eq.~\ref{MLAPLACE} for $\tilde{T}(s)$, explicitly establishing the
mean-field equivalence between the Bellman-Harris approach and our
kinetic theory.  Note that in the Bellman-Harris formulation, the
waiting-time distributions of either fission or death have the same
distribution $g(a)$. In our kinetic theory, these rates can have
distinct distributions, $\beta_n(a)$ and $\mu_n(a)$, and can also
depend on population size, providing much greater flexibility.


\section{A Cell Division Model}
\label{CELLSECTION}

We now consider explicit calculations for a simple fission-only model
of cell division where cell cycle times are rescaled to be
$\Gamma$-distributed with unit mean and variance
$\frac{1}{\alpha}$. This $\Gamma$-distribution and its Laplace
transform $\tilde{g}(s)$ are explicitly

\begin{equation}
g(t) = \frac{\alpha^\alpha}{\Gamma(\alpha)}t^{\alpha-1}e^{-\alpha t},
\hspace{10mm}
\tilde{g}(s) = \left (\frac{\alpha}{\alpha+s}\right )^\alpha.
\end{equation} 
Equation \ref{PHILAPLACE} for $B(t)$ can then be solved to yield

\begin{equation}
B(t) = \mathcal{L}^{-1}_t\left(\frac{\alpha^\alpha}{(s+\alpha)^\alpha-2\alpha^\alpha}\right)
=\alpha e^{-\alpha t}\mathcal{L}^{-1}_{(\alpha t)}\left(\frac{1}{s^\alpha-2}\right).
\label{PHILAPINV}
\end{equation}
The inverse Laplace transform is detailed in Appendix B and 
involves contour integration that yields

\begin{equation}
B(t) = -\frac{\alpha}{\pi}\int_0^{\infty}\frac{e^{-\alpha t(r+1)}
r^\alpha\sin(\pi \alpha)}{r^{2\alpha}-4r^\alpha\cos(\pi \alpha)+4}\dd r
+\sum_{n=-\left \lfloor\frac{\alpha}{2}\right\rfloor}^{\left \lfloor\frac{\alpha}{2}\right\rfloor}
2^{\frac{1}{\alpha}-1}e^{{(2^{\frac{1}{\alpha}}\cos(\frac{2n\pi}{\alpha})-1)}\alpha t}
\cos\left(2^{\frac{1}{\alpha}}\alpha t\sin\left(\frac{2n\pi}{\alpha}\right)+\frac{2n\pi}{\alpha}\right).
\label{PHIFORMULA}
\end{equation}

Similarly, from Eq.~\ref{MLAPLACE} we have

\begin{equation}
T(t) = \mathcal{L}_t^{-1}\left({\frac{1}{s}\frac{(s+\alpha)^\alpha-\alpha^\alpha}
{(s+\alpha)^\alpha-2\alpha^\alpha}}\right)=e^{-\alpha t}\mathcal{L}_{(\alpha t)}^{-1}
\left({\frac{1}{s-1}\frac{s^\alpha-1}{s^\alpha-2}}\right),
\label{BROM2}
\end{equation}
which can also be evaluated via a similar Bromwich integral:

\begin{align}
T(t) = & \displaystyle \frac{1}{\pi}\int_0^{\infty}\frac{e^{-\alpha t(r+1)}}{r+1}
\frac{r^\alpha\sin(\pi \alpha)}{r^{2\alpha}-4r^\alpha\cos(\pi \alpha)+4}\dd r \nonumber\\
 & +\sum_{n=-\left \lfloor\frac{\alpha}{2}\right\rfloor}^{\left \lfloor\frac{\alpha}{2}\right\rfloor}
\frac{2^{\frac{1}{\alpha}}}{2\alpha}e^{{(2^{\frac{1}{a}}\cos(\frac{2n\pi}{\alpha})-1)}\alpha t}
\frac{2^{\frac{1}{\alpha}}\cos(2^{\frac{1}{\alpha}}
\sin(\frac{2n\pi}{\alpha})\alpha t)-\cos(2^{\frac{1}{\alpha}}
\sin(\frac{2n\pi}{\alpha})\alpha t+\frac{2n\pi}{\alpha})}
{2^{\frac{2}{\alpha}}-2^{1+\frac{1}{\alpha}}\cos(\frac{2n\pi}{\alpha})+1}.
\label{MFORMULA1}
\end{align}

\begin{figure}[t!]
\begin{center}
\includegraphics[width=15cm]{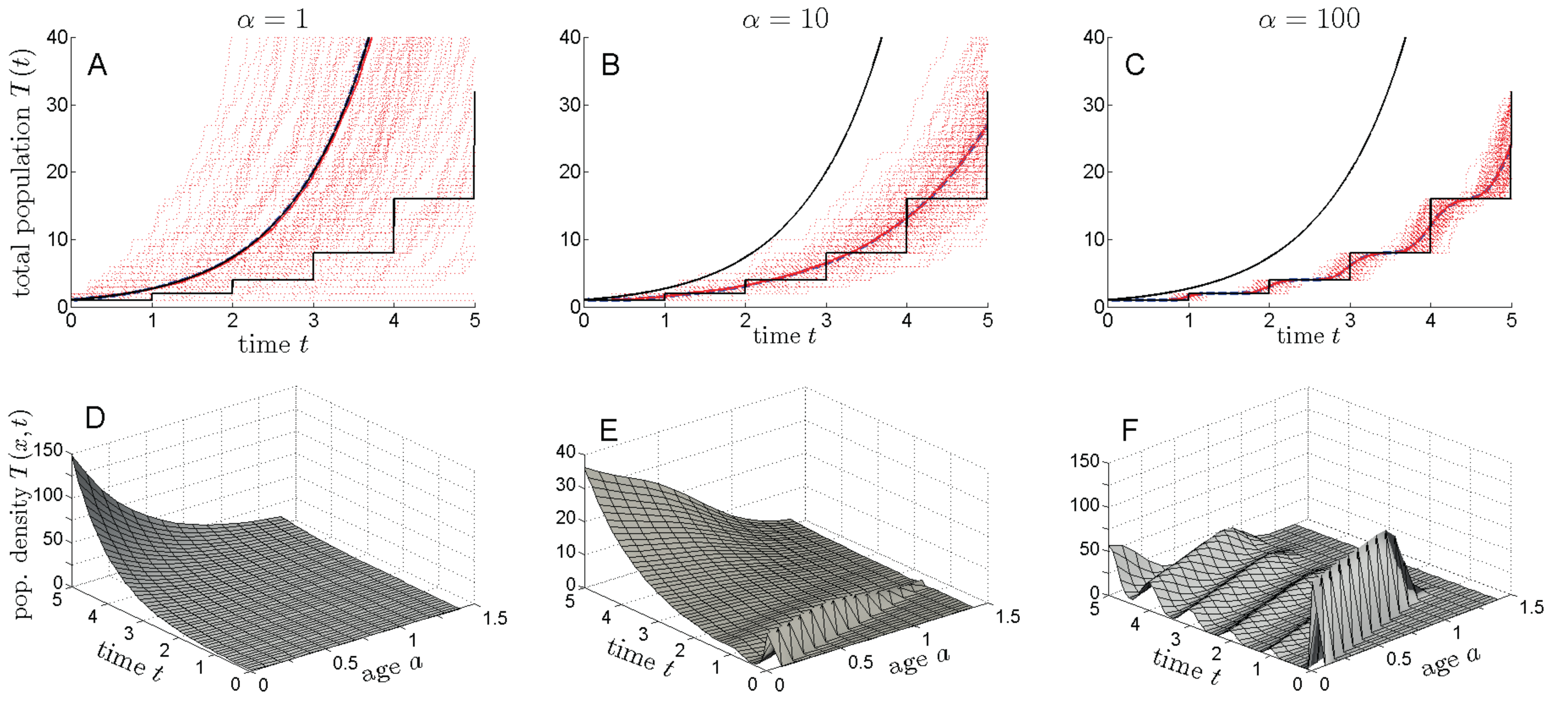}
\caption{Mean population and age-time distributions for a fission-only
  process with $\Gamma$-distributed branching times.  \textsf{A},
  \textsf{B}, and \textsf{C} show mean populations as a function of
  time for dispersion values $\alpha=1$, $\alpha=10$, and $\alpha =
  100$, respectively.  Red dotted trajectories are realizations of
  simulations, while the solid red line is the mean.  The blue dashed
  curve is the mean population $T(t)$ computed from Eq.~\ref{MFORMULA1} and
  is nearly indistinguishable from the red solid curve. The upper and
  lower black lines correspond to the continuous-time Markovian
  fission process and the discrete-time Galton-Watson process,
  respectively.  \textsf{D},\textsf{E}, and \textsf{F} depict the
  corresponding mean age-distributions $T(x,t)$ computed from
  Eq.~\ref{TFULLFORMULA} but plotted as functions of time $t$ and age
  $a$.\label{YULEPLOT}}
\end{center}
\end{figure}

For $\alpha = 1$, $g(t) = e^{-t}$ is exponentially distributed, and we
find the simple growth law $T(t) = e^{t}$, which is equivalent to the
result $\mathrm{E}(Y_{[0,\infty]}) = e^{\beta t}$ found earlier in
Subsection \ref{MOMENTEQNS}. As $\alpha \to \infty$, the
gamma-distribution sharpens about unity and the process becomes more
discrete-like in time.  Figs.~\ref{YULEPLOT}\textsf{A,B,C} show that
as $\alpha$ is increased, the mean population size $T(t)$ tends
towards that given by the discrete-time Galton-Watson step process, as
would be expected.  In Figs.~\ref{YULEPLOT}\textsf{D,E,F}, we have
used the expression for $B(t)$ in Eqs.~\ref{TFULLFORMULA} and
\ref{PHIFORMULA} to give the mean age-time distribution $T(x,t)$. Note
that unlike the solution to the Bellman-Harris equation shown in
Figs.~\ref{YULEPLOT}\textsf{A,B,C}, the mean density $T(x,t)$
(Eq.~\ref{TFULLFORMULA}) resolves age structure. 

%


\section{Spatial Models}

We now illustrate how our kinetic (in age) theory can be generalized
to include spatial motion such as diffusion and convection.  We will
follow the approaches described in Webb \cite{WEBB2008} for
incorporating spatial effects in age-structured simple birth-death
processes. Since these methods are adaptations of the McKendrick-von
Foerster equation, they are deterministic and ignore stochastic
fluctuations in population size. In a manner similar to how the
McKendrick-von Foerster equation was extended to the stochastic domain
using Eq.~\ref{RHO0}, in this section we briefly outline how to
generalize the age-structured spatial process discussed in
\cite{WEBB2008} to incorporate stochasticity.

Consider a simple budding-mode birth-death process such that
$\hat{\rho}_n(\b_n;\q_n;t)$ is the probability density for a
population containing $n$ randomly labelled individuals with TOBs
$\b_n$ \textit{and} positions $\q_n$. Although
$\hat{\rho}_n(\b_n;\q_n;t)$ is again invariant under permutations of
the elements, the relative orders of $\b_n$ and $\q_n$ must be
preserved: $\hat{\rho}_2(b_1,b_2;q_1,q_2;t) = \hat{\rho}_2
(b_2,b_1;q_2,q_1;t)$, for example. The TOB arguments in this density
can be readily transformed to age (rather than TOB) by a
transformation of the type given in Eq.~\ref{AANDB}. For ease of
presentation, we assume a one-dimensional system; generalizations to
higher spatial dimensions are straightforward. We further suppose that
individuals are undergoing identical, independent diffusion processes
with diffusion constant $D$. Examples of other spatial processes that
may be incorporated can be found in \cite{WEBB2008}. We suppose that
$\beta_n(a;q)$ and $\mu_n(a;q)$ are birth and death rates for
any individual with age $a$ and at spatial position $q$ in a population
of size $n$. Finally, the initial position of each newborn is
determined by the position of the parent at the time of birth.  The
extended theory is described by the following kinetic equation for
$\hat{\rho}_n(\b_n;\q_n;t)$:

\begin{align}
\frac{\partial \hat{\rho}_{n}(\b_n;\q_n;t)}{\partial t} = & -\hat{\rho}_{n}(\b_n;\q_n;t) \sum_{i=1}^n\gamma_{n}(t-b_i,q_i)+ D\sum_{i=1}^n\frac{\partial^2}{\partial q_i^2}\hat{\rho}_n(\b_n;\q_n;t)
\nonumber\\
&+(n+1)\int_{-\infty}^{t}\dd y 
\int_{\mathbb{R}} \dd q' \hspace{1mm}\hat{\rho}_{n+1}(\b_n,y;\q_n,q';t)\mu_{n+1}(t-y,z).
\label{SPACEEQ}
\end{align}
%
%
The corresponding boundary condition capturing the influx of newborn
individuals is

\begin{equation}
\rho_n(\b_
{n-1},t;\q_n;t) =
\frac{1}{n}\sum_{i=1}^{n-1}\rho_{n-1}(\b_{n-1};\q_{n-1};t)\beta(t-b_i,q_i)\delta(q_n-q_i),
\end{equation}
which differs slightly from that in Eq.~\ref{BCRHO}. In the original
formulation, we do not track which individual is the parent of a
newborn, whereas here the newborn has the same position ($q_n$) as the
parent ($q_i$), setting its identity as the $i^\mathrm{th}$
individual. In addition to a boundary condition, Eq.~\ref{SPACEEQ}
requires an initial condition $\rho_n(\b_n;\q_n;0)$ to specify both
the initial TOB and initial position of individuals.

As with our earlier analyses, we first express $\rho_n$ in terms of
$\rho_{n+1}$ by introducing the propagator
$U_n(\b_n;\q_{n};t_0,t)=\exp\left[-\sum_{i=1}^n\int_{t_0}^t
\gamma_n(s-b_i,q_i)\dd s\right]$, which enables us to transform
Eq.~\ref{SPACEEQ} to an inhomogeneous heat equation for the function
$U_n^{-1}\rho_n$

\begin{equation}
\frac{\partial}{\partial t}\left[U_n^{-1}(\b_n;\q_{n};t_0,t)\rho_n\right]=
D\sum_{j=1}^n\frac{\partial^2}{\partial
  q_j^2}\left[U_n^{-1}\rho_n\right] +
(n+1)U_n^{-1}\int_{-\infty}^{t}\hspace{-4mm}\dd y \int_\mathbb{R}
\!\!\! \dd z\hspace{1mm}\rho_{n+1}(\b_n,y;\q_n,z;t)\mu_{n+1}(t-y,z),
\end{equation}
whose solution can be expressed as \cite{HEATEQN}

\begin{align}
\rho_n(\b_n;\q_n;t) = & U_n(\b_n;\q_{n};t_0,t)\int_{{\mathbb R}^n}\dd \q'_n 
\HS N_{\q_n}(\q'_n,2D(t-t_0)I_n)\rho_n(\b_n;\q'_n;t_0) 
\nonumber\\
& +(n+1)\int_{t_0}^t \dd s \HS U_n(\b_n;\q_{n};s,t)\int_{{\mathbb R}^m}
\dd \q'_n N_{\q_n}(\q'_n,2D((t-t_0)-s)I_n)
\nonumber\\
& \hspace{30mm}\times\int_{-\infty}^s \dd y \int_{\mathbb{R}} \dd z 
\,\rho_{n+1}(\b_n,y;\q'_n;z;s)\mu_{n+1}(s-y,z).
\label{RECCSPACE}
\end{align}
Here, $I_n$ denotes the $n\times n$ identity matrix and
$N_{\q}(\x,\Sigma)$ is the multivariate normal density for the vector
$\q$ arising from a distribution with mean $\x$ and covariance
$\Sigma$. This result expresses $\rho_{n}$ in terms of $\rho_{n+1}$
and is analogous to Eq.~\ref{RECC}. This solution is valid
provided $t_0>\max\{\x\}$; for $t_0=\max\{\x\}$, we must invoke the
boundary condition. One can then use Eq.~\ref{RECCSPACE} and the
boundary condition to search for explicit solutions in much the same
way as we did for our spatially independent kinetic theory.
%
In the companion paper, we derive the mean-field equations for this
spatial kinetic theory using quantum field theoretic methods developed
by Doi and Peliti \cite{Greenman2}.


\section{Summary and Conclusions}

We have developed a complete kinetic theory for age-structured
birth-death and fission-death processes that allow for systematic and
and self-consistent incorporation of interactions at the population
level. Our overall result in \cite{Greenman1} which we extend here is
the derivation of a kinetic theory for stochastic age-structured
populations.  The kinetic equations can be written in terms of a
BBGKY-like hierarchy (or a double hierarchy in the case of
fission). Methods of approximation and closure typically employed in
gas/liquid kinetic theory, plasma physics, or fluid dynamics can then
be applied.

The analysis presented in this paper provides three additional
specific mathematical results. Firstly, in Eq.~\ref{MVFG}, we have
shown that the factorial moments of the age structure can be described
by an equation that naturally generalizes the McKendrick-von Foerster
equation.  In particular, for population-independent birth, death, and
fission rates we can determine the variance of the population size for
specific age groups in a population, something that was not previously
feasible without some form of approximation.

Secondly, in Eqs.~\ref{RHO0} and \ref{BCRHO}, we develop a complete
probabilistic description of a population undergoing a binary fission
and death process. Although a general analytic solution to these
systems can be written down (Eq.~\ref{GENSOL}), it is difficult to
calculate and further work is needed to identify analytic techniques
or numerical schemes that can more readily provide solutions. The
methods we have introduced can also be viewed as a continuum limit of
matrix population models.

Thirdly, we also outlined how to incorporate spatial dependence of
birth and death into our age-structured kinetic theory. We considered
only the simplest model of free diffusion in which individuals to not
interact spatially. Spatially-mediated interactions can be
incorporated by way of a ``collision operator'' in a full theory that
treats both age and space kinetically.

All of our results can also be derived using from
quantum field theoretical approaches \cite{Doi1,Doi2,Peliti}, which
are described in detail in a parallel paper \cite{Greenman2}. Such
methods provide alternative machinery to analyze the statistics of
age- and space-structured populations and may provide new avenues for
calculation.

Finally, we note that the overall structure of our model is
semi-Markov.  That is, birth, death, and fission rates depend on
only the time since birth of an individual and not on, for example,
the number of generations removed from a founder.  Such lineage aging
processes are often important in cell biology (\textit{e.g.}, the
Hayflick limit \cite{HAYFLICK1965}) and would require extension of our
state space to include generational class \cite{ZILMAN2010}. These
extensions will be explored in future work.
 

\section*{Appendix A: Second Factorial Moment Derivation}

We outline how to derive Eq.~\ref{APPASTUFF}.  Assume the initial
population is described by $X^{(1)}(a;0)=\lambda e^{-\lambda a}$ and
$X^{(2)}(a,b;0) = 0$.  Note that $X^{(1)}$ is just the solution to the
McKendrick-von Foerster equation given by the expression in
Eq.~\ref{MVFSOL}.  We can determine $X^{(2)}$ via Eq.~\ref{FMOMSOL} if
we are able to identify the boundary condition $B(a,t) \equiv
X^{(2)}(a,0;t) \equiv X^{(2)}(0,a;t)$. After setting $m=\min\{a,b,t\}$
in Eq.~\ref{FMOMSOL}, we substitute the expressions for $X^{(2)}$ into
the boundary condition Eq.~\ref{FMBC} to give the following equation
for $B(a,t)$:

\begin{equation}
B(a,t) = \frac{\beta}{2}X^{(1)}(a;t)+
\beta\begin{cases}
\int_0^tB(a-b,t-b) \dd b, & t<a,\\
\int_0^aB(a-b,t-b)\dd b
+\int_a^\infty B(b-a,t-a)\dd b, & t>a.
\end{cases}
\label{BCAPP}
\end{equation}

An expression for $B(a,t)$ in the region $t<a$ can be obtained by
solving along characteristics such as those portrayed in
Fig.~\ref{TIMEAGE}. We first define $C(\alpha,\tau)=B(a,t)$, where
$\alpha=a-t, \tau = t$, so that

\begin{equation}
C(\alpha,\tau) = \frac{\beta}{2}X^{(1)}(\alpha+\tau;\tau)+
\beta\int_0^t C(\alpha,\tau-b)\dd b.
\end{equation}
A Laplace transform with respect to $\tau$ can then be used to
find $B(a,t) = \frac{\beta\lambda}{2} e^{-\lambda  a}e^{(\lambda+\beta)t}$.

For $t>a$, note that the second integral in Eq.~\ref{BCAPP} extends
into the region $t<a$, for which we now have an expression. Upon
separating the integral into two parts, and similarly defining
$C(\alpha,\tau)=B(a,t)$, where $\alpha=a, \tau = t-a$ along
characteristics, we find

\begin{equation}
C(\alpha,\tau) = \frac{\beta^2}{2}e^{\beta\tau}
+\beta\int_0^\alpha C(b,\tau)\dd b
+\beta\int_0^\tau C(b,\tau-b)\dd b
+\frac{\beta\lambda}{2}\int_{\tau+\alpha}^\infty
e^{-\lambda(b-\alpha)}e^{(\lambda+\beta)\tau}\dd b.
\end{equation}
A double Laplace transform in variables $\alpha$ and $\tau$ results in:

\begin{equation}
\hat{C}(u,v) = \frac{\beta}{u}\left(\hat{C}(u,v)
+\hat{C}(v,v)\right)+\frac{\beta^2}{u}\frac{1}{v-\beta},
\end{equation}
from which we find $\hat{C}(v,v) =
\frac{\beta^2}{(v-\beta)(v-2\beta)}$ and so $\hat{C}(u,v) =
\frac{\beta^2}{(u-\beta)(v-2\beta)}$. A double Laplace inversion then
gives $B(a,t) = \beta^2 e^{-\beta a}e^{2\beta t}$, from which $X^{(2)}$ can
be uniquely determined from Eq.~\ref{FMOMSOL}.


\section*{Appendix B: Bromwich Integral Calculation}

\begin{figure}[t!]
\begin{center}
\includegraphics[width=6cm]{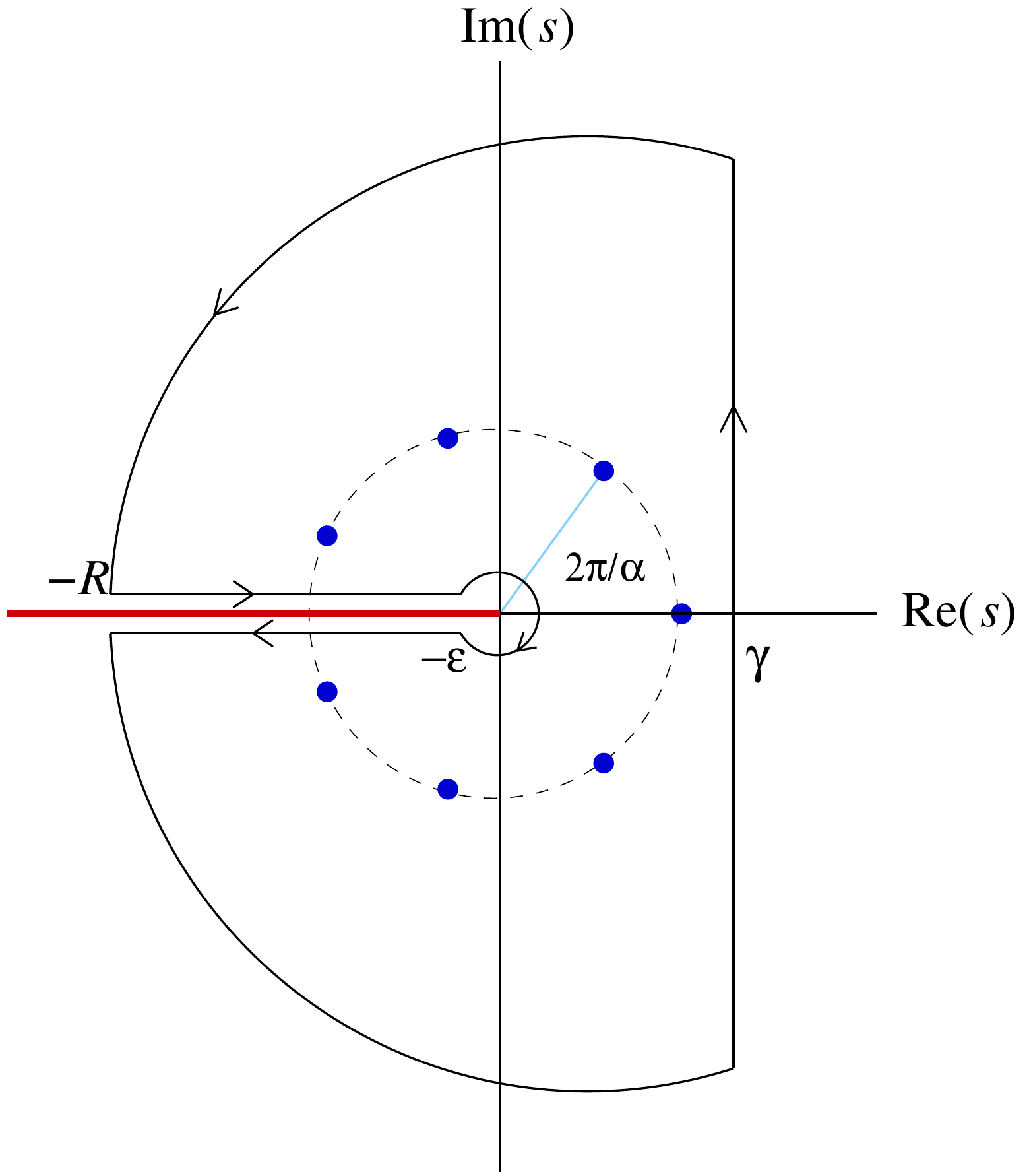}
\caption{Bromwich integral for calculating the inverse Laplace
  transform in Eq.~\ref{BROM1}. The integral along $\gamma$ is
  evaluated using the residues at the poles and 
the integrals along the branch cut in Cauchy's theorem.\label{CONTOUR}}
\end{center}
\end{figure}

Since the inverse Laplace transform provided by the Bromwich integral

\begin{equation}
\mathcal{L}^{-1}_{t}\left(\frac{1}{s^\alpha-2}\right)=\frac{1}{2\pi i}
\int_{\gamma-i\infty}^{\gamma+i\infty}\frac{e^{st}}{s^\alpha-2}\dd s
\label{BROM1}
\end{equation}
involves a branch point at $s=0$, we construct a branch cut along the
negative real axis and define $s=re^{i\theta}$ where $\theta \in
(-\pi,\pi)$. The denominator $s^{\alpha}-2$ also produces poles at
$s=2^{\frac{1}{\alpha}}e^{i\frac{2n}{\alpha}}$ where $n$ is an integer
with $|n| \le \left \lfloor{\frac{\alpha}{2}}\right \rfloor$. The
contour required for the Bromwich integral is shown in
Fig.~\ref{CONTOUR} and is evaluated using Cauchy's residue
theorem. 

The integrals around the outer perimeter and the origin contribute
zero in the limit as $R \to \infty$ and $\varepsilon \to 0$. The
branch cuts and poles provide the nonzero contributions.  First,
consider the integrals along the branch cut. Writing the variable $s$
as $re^{i\theta}$, for $\theta = \pm \pi$, we integrate $\frac{1}{2\pi
  i}\frac{e^{st}}{s^{\alpha}-2}$ along the two sides to give

\begin{equation}
\frac{1}{2\pi i}\int_\infty^0\frac{e^{-rt}(\dd r e^{i\pi})}{r^\alpha e^{i\pi\alpha}-2}
+\frac{1}{2\pi i}\int_0^\infty\frac{e^{-rt}(\dd r e^{-i\pi})}{r^\alpha e^{-i\pi\alpha}-2}
= -\frac{1}{\pi}\int_0^\infty\frac{e^{-rt}r^\alpha\sin(\pi\alpha)\HS\dd r}{r^{2\alpha}-4r^\alpha\cos(\pi\alpha)+4}.
\end{equation}
Next, we need to consider the poles at positions
$s=2^\frac{1}{\alpha}e^\frac{2n\pi i}{\alpha}$ for $|n| \le
\lfloor{\frac{\alpha}{2}}\rfloor$. L'H\^{o}pital's rule leads to

\begin{equation}
\lim_{s \to 2^\frac{1}{\alpha}e^\frac{2n\pi i}{\alpha}}
\left\{\frac{s-2^\frac{1}{\alpha}e^\frac{2n\pi i}{\alpha}}{s^\alpha-2}\right\} = \lim_{s \to 2^\frac{1}{\alpha}e^\frac{2n\pi i}{\alpha}}\left\{\frac{1}{\alpha s^{\alpha-1}}\right\}
= \alpha^{-1}2^{\frac{1}{\alpha}-1}e^{\frac{2n\pi i}{\alpha}}.
\end{equation}
If $r_n$ is the residue for the function $\frac{e^{st}}{s^\alpha-2}$
at the pole $s=2^\frac{1}{\alpha}e^\frac{2n\pi i}{\alpha}$, we can
write

\begin{equation}
r_n+r_{-n} = 2\mathrm{Re}\left\{\alpha^{-1}2^{\frac{1}{\alpha}-1}e^{\frac{2n\pi i}{\alpha}}e^{2^{\frac{1}{\alpha}}e^{\frac{2n\pi i}{\alpha}}t}\right\} = \frac{2^{\frac{1}{\alpha}}}{\alpha}e^{2^{\frac{1}{\alpha}}\cos\left(\frac{2n\pi}{\alpha}t\right)}
\cos\left(2^{\frac{1}{\alpha}}\sin\left(\frac{2n\pi}{\alpha}\right)+\frac{2n\pi}{\alpha}\right).
\end{equation}
Combining the contributions from the branch cut and the residues
results in $\mathcal{L}_{(t)}^{-1}\left(\frac{1}{s^\alpha-2}\right)$,
which, when substituted into Eq.~\ref{PHILAPINV}, gives the final
result in Eq.~\ref{PHIFORMULA}.

The derivation for the Laplace inversion in Eq. \ref{BROM2} is
similar. Note that the value $s=1$ is a removable singularity and the
same set of poles and paths for branch cut integrals needs to be
considered. The details are left to the reader.


\section*{Acknowledgements}

This research was supported in part by the National Science Foundation
under Grant No. NSF PHY11-25915 to KITP and by the Gordon and Betty
Moore Foundation under Award No. 2919 to the KITP. TC is also
supported by the US National Institutes of Health through grant R56
HL126544, the NSF through grant DMS-1516675, and the Army Research
Office through grant W911NF-14-1-0472.

%
%


\bibliographystyle{spmpsci}      

\bibliography{refs_Fission}

%
%

\end{document}